\documentclass[12pt,a4paper]{amsart}
\usepackage{amsmath,amsfonts,amssymb,amsthm,amscd}
\usepackage{dsfont,mathrsfs,enumerate}
\usepackage{curves,epic,eepic}

\usepackage[pdftex]{graphicx}
\usepackage{psfrag}
\usepackage{epic}
\usepackage{eepic}
\usepackage[matrix,arrow,curve]{xy}
\input{epsf}
\usepackage{tikz}
\usepackage{xcolor}
\usetikzlibrary{calc}
\usetikzlibrary{intersections}
\usetikzlibrary{through}
\usetikzlibrary{backgrounds}

\newtheorem{thm}{Theorem}[section]

\newtheorem{prop}[thm]{Proposition}
\newtheorem{coro}[thm]{Corollary}
\theoremstyle{definition}
\newtheorem{defi}[thm]{Definition}
\newtheorem{rem}[thm]{Remark}
\newtheorem{ex}[thm]{Example}

\def\T{\mathds T}
\def\F{\mathds F}

\def\Z{\mathds Z}
\def\Q{\mathds Q}
\def\R{\mathds R}
\def\C{\mathds C}
\def\phi{\varphi}

\def\P{\mathds P}

\def\<{{\langle}}
\def\>{{\rangle}}

\DeclareMathOperator{\Hom}{Hom}

\DeclareMathOperator{\ord}{Min}

\begin{document}

\title[Canonical models of toric hypersurfaces]{
Canonical models of toric hypersurfaces}

\author[Victor V. Batyrev]{Victor V. Batyrev}
\address{Mathematisches Institut, Universit\"at T\"ubingen,
Auf der Morgenstelle 10, 72076 T\"ubingen, Germany}
\email{batyrev@math.uni-tuebingen.de}

\begin{abstract}
Let $Z$ be a 
nondegenerate hypersurface in
$d$-dimensional torus $(\mathbb{C}^*)^d$ defined 
by a Laurent polynomial $f$ with a $d$-dimensional 
Newton polytope $P$. The subset $F(P) \subset P$ consisting 
of all points in $P$ having integral distance at least $1$ to  all integral supporting hyperplanes of $P$ is called the 
Fine interior of $P$. If $F(P) \neq \emptyset$  
we construct    
a unique projective    
model $\widetilde{Z}$  of $Z$ having at worst canonical singularities  and   obtain  
minimal models $\widehat{Z}$  of $Z$ by  
crepant morphisms  $\widehat{Z} \to  
\widetilde{Z}$. 
We show that the Kodaira dimension $\kappa =\kappa(\widetilde{Z})$ equals 
$\min \{ d-1, \dim F(P) \}$ and the general  fibers in  the    
Iitaka fibration of the canonical model  $\widetilde{Z}$ are
non\-degenerate $(d-1-\kappa)$-dimensional
toric hypersurfaces of  Kodaira dimension $0$. Using 
$F(P)$, we obtain  a simple combinatorial 
formula for the intersection number $(K_{\widetilde{Z}})^{d-1}$. 
\end{abstract}


\maketitle

\thispagestyle{empty}
\section{Historical motivations}

Let $M \cong \Z^d$ be a $d$-dimensional lattice that  we identify with the character group of  a $d$-dimensional algebraic torus 
$\T^d \cong (\C^*)^d$ whose affine coordinate ring $\C[\T^d] \cong \C[ t_1^{\pm 1}, \cdots, t_d^{\pm 1} ]$ consists 
of Laurent polynomials 
\[ f({\bf t}) = \sum_{ {\bf m} \in M} c_{\bf m} 
{\bf t}^{\bf m}, \; \, c_{\bf m} \in \C, \]
where $c_{\bf m} = 0$ for all but finitely many ${\bf m} \in M$. The Newton polytope $P:=Newt(f)$ of $f({\bf t})  \in \C[\T^d]$ is the convex hull of all ${\bf m} \in M$ such that $c_{\bf m} \neq 0$. We assume that $\dim P =d$ and  write  
$ f({\bf t}) =\sum_{ {\bf m} \in A} c_{\bf m} 
{\bf t}^{\bf m}$, where $A \subset M$ is a finite subset with $P = 
{\rm Conv}(A)$. Following  Khovanski\v{i} \cite{Kho77},  a  Laurent polynomial $f \in \C[\T^d]$ and  an  
affine hypersurface $Z:= \{ f({\bf t} )= 0 \} \subset \T^d$ are 
called {\em nondegenerate} if the  
Zariski closure $\overline{Z}_P$ of $Z$ in  projective toric variety $V_P$ defined by the normal fan $\Sigma_P$ has transversal intersections with all $k$-dimensional $\T^d$-orbits 
$\T_Q \subset V_P$ corresponding to  $k$-dimensional
 faces $Q \preceq P$ $(k \geq 1)$.

The set of coefficients $\{ c_{\bf m} \}_{{\bf m } \in A} 
\in \C^{|A|}$ defining {\em nondegenerate Laurent polynomials} 
$ f({\bf t}) =\sum_{ {\bf m} \in A} c_{\bf m} 
{\bf t}^{\bf m}$ is a Zariski open dense subset $U_A \subset \C^{|A|}$ that can be explicitly described by the nonvanishing 
condition for the  
{\em principal $A$-determinant} $E_A$ introduced  
by Gelfand, Kapranov and Zelevinski \cite{GKZ94}. 

According to Khovanski\v{i} \cite{Kho77,Kho78,Kho83}, one can resolve 
singularities of the projective hypersurface $\overline{Z}_P \subset V_P$ and 
obtain a smooth projective model $\overline{Z}$ of  a  nondegenerate  
hypersurface $Z \subset \T^d$ as
 Zariski closure $\overline{Z}$  of $Z$ in 
some  smooth projective equivariant torus embedding $\T^d \hookrightarrow V$ corresponding to  
a regular simplicial refinement $\Sigma$ 
of the normal fan $\Sigma_P$. 
Obviously, there exist countably many possibilities to choose such a refinement $\Sigma$.

The Minimal Model Program (MMP), or Mori theory,  is aimed at constructing  the most "economic" projective representatives in each birational class of algebraic varieties. 
Such representatives are called 
{\em minimal models} and they are main objectives   
of  MMP \cite{KM98,Mat02,Bir18a}.  

The present paper is inspired by the author's construction of 
$(d-1)$-dimensional minimal Calabi-Yau varieties $\widehat{Z}$ 
using projective  compactifications  of 
 nondegenerate affine hypersurfaces $Z \subset \T^d$ defined by Laurent polynomials whose Newton polytopes $P$ are {\em reflexive} \cite{Bat94}. These minimal Calabi-Yau varieties allow to 
 produce many topologically 
different examples of smooth Calabi-Yau $3$-folds, 
compute their 
Hodge numbers by combinatorial formulas and perfectly illustrate the mirror symmetry phenomenon in theoretical physics \cite{CK99, KS98,KS02,AGHJN15}. We remark that the normal fan $\Sigma_P$ of any reflexive 
polytope $P$ defines a Gorenstein toric Fano variety $V_P$. 
One can choose an  "economic" simplicial refinement $\widehat{\Sigma}$ of $\Sigma_P$ such that the corresponding toric varie\-ty $\widehat{V}$ has  at worst 
terminal singularities and the corresponding 
proper birational toric morphism $\widehat{V} \to V_P$ is {\em crepant}. Such a refinement $\widehat{\Sigma}$ may be not unique. The Zariski closure $\widetilde{Z} :=\overline{Z}_P$  of the 
affine hypersurface $Z \subset \T^d$ in $\widetilde{V} := V_P$ is a Calabi-Yau 
hypersurface $\widetilde{Z}$ with at worst Gorenstein canonical singularities. The minimal Calabi-Yau models  $\widehat{Z}$  of $Z$ are partial crepant desingularizations of $\widetilde{Z}$ obtained as Zariski closures of affine hypersurfaces $Z \subset \T^d$ in the 
simplicial toric varieties $\widehat{V}$.

A generalization of the 
above construction to minimal models of  arbitrary  
nondegenerate toric hypersurfaces $Z \subset \T^d$ has been 
suggested by Shihoko Ishii  \cite{Ish99},   who  
showed  that  minimal models of $Z$  can be {\em always} obtained as Zariski closures $\widehat{Z}$
of $Z$ in some appropriately chosen 
simplicial toric varieties $\widehat{V}$ having at worst terminal singularities.  Moreover, Ishii proposed two ways for constructing such a  toric variety $\widehat{V}$: the traditional method and a shorter one.  

The traditional method of MMP begins   
from a smooth projective 
model  $\overline{Z} \subset V$ obtained in some smooth toric variety $V$ by Khovanski\v{i}'s method. 
The traditional MMP-method suggests a sequence  of birational   
modifications of the pair $(V, \overline{Z})$ using the combinatorial framework of toric Mori theory developed by Miles Reid  \cite{Rei83} (see also \cite{Fuj03,FS04} and \cite[\S 15]{CLS11}). It is useful to keep in mind the short exact
sequence describing the canonical sheaf ${\mathcal K}_{\overline{Z}}$  on $\overline{Z}$ by the adjunction formula:
\[ 0 \to {\mathcal O}_V(K_V ) \to {\mathcal O}_V(K_V + \overline{Z}) \to {\mathcal K}_{\overline{Z}} \to 0. \]
We remark  that  the Kodaira dimension $\kappa := 
\kappa(\overline{Z})$ of 
the smooth hypersurface $\overline{Z} \subset V$ is nonnegative if and only if the linear systems corresponding 
to  some multiples of the  adjoint divisor $K_V + \overline{Z}$ on $V$ are not empty, i.e., that the adjoint class 
 $[K_V + \overline{Z}]$ belongs to the closed cone $C_{\rm eff}(V) \subset H^2(V, \R)$ of effective 
divisors of $V$. 

If the adjoint divisor 
$K_V + \overline{Z}$ is nef, then
$\overline{Z}$ is a minimal model. Otherwise there exists 
a toric extremal ray $R=\R_{\geq 0}v \in H_2(V,\R)$ 
with the negative intersection number $(v, K_V+ \overline{Z}) <0$. Applying  an elementary  
 birational toric modification  from the toric Mori theory, one can obtain   
another  pair $(V', \overline{Z}')$ consisting of a simplicial toric variety $V'$ having at worst terminal 
singularities and another projective model  $\overline{Z}' \subset V'$ of $Z$ etc. Eventually, after finitely many elementary  birational  
toric modifications,   
one  comes to a projective simplicial toric variety 
$\widehat{V}$ with at worst 
terminal singularities containing a 
$\Q$-Cartier divisor $\widehat{Z}$ 
such that the adjoint divisor 
$K_{\widehat{V}} +  \widehat{Z}$ is nef.  
Then the projective hypersurface 
 $\widehat{Z} \subset \widehat{V}$ is a minimal model of $Z$. 
Unfortunately, the traditional MMP-method does not allow to say much about the resulting simplicial toric variety $\widehat{V}$ and 
about the minimal model $\widehat{Z}$ ifself. 

For this reason,  Ishii suggested a more informative and short way for constructing  
the  simplicial 
toric variety $\widehat{V}$. It uses  
the rational polytope $F(P) \subset M_\R$  correspon\-ding to the adjoint divisor $K_V+ \overline{Z}$ on some smooth 
projective toric 
variety $V$ obtained by Khovanski\v{i}'s method using a refinement 
$\Sigma$ of the normal fan of $P$.\footnote{Shihoko Ishii denotes the rational polytope $F(P)$ by 
$\square_h$ (see \cite[3.3]{Ish99}).} It is easy to see that the rational polytope 
$F(P)$  is independent of  such a refinement $\Sigma$.  We note that the  polytope $F(P)$ was earlier 
introduced in the doctoral thesis of Jonathan 
Fine \cite{Fin83} and, 
following Reid \cite[Appendix to \S 4]{Rei87}, 
we call the polytope $F(P)$ the {\em Fine interior} of $P$.\footnote{Jonathan Fine uses  for the rational polytope 
$F(P)$ the name 
{\em heart of $P$} \cite[\S 4]{Fin83}.}  An 
important observation of Ishii  is that all $1$-dimensional 
cones in  the simplicial fan $\widehat{\Sigma}$ defining  
the toric variety $\widehat{V}$ are uniquely determined by  $P$. They are spanned  by the lattice vectors 
from a special finite  set 
$$S_F(P):= \{\nu_1, \ldots, \nu_p \} = \widehat{\Sigma}[1] \subset \Sigma[1]$$ 
characterized by the property 
\[ \nu \in S_F(P) \Leftrightarrow \min_{x \in F(P)}\langle x, \nu \rangle  = 
\min_{x \in P}\langle x, \nu \rangle +1. \] 
 Ishii called  the lattice vectors $\nu \in S_F(P)$ 
{\em contributing to the Fine interior} $F(P)$. 
In the present  paper,  we call  the finite set $S_F(P)$  the 
{\em support of the Fine interior} of $P$.  

Ishii suggested   to 
construct  the  simplicial fan $\widehat{\Sigma}$ with 
$\widehat{\Sigma}[1] = S_F(P)$ as normal fan 
of some simple $d$-dimensional 
polytope $\square(\varepsilon)$ whose primitive inward pointing facet normals are exactly the lattice vectors 
$\nu_1, \ldots, \nu_p \in S_F(P)$. Note that 
there exists a natural 
bijection between all projective simplicial fans $\widehat{\Sigma}$ 
with $\widehat{\Sigma}[1] = S_F(P)$ 
and maximal $(p-d)$-dimensional $GKZ$-cones, or Mori chambers,
in the moving cone  of $\widehat{V}$ (see  \cite{OP91}, \cite[Appendix A]{HKP06}, \cite[\S12]{CLS11}).
In order to obtain a minimal model $\widehat{Z} \subset \widehat{V}$ one has to choose the simple polytope $\square(\varepsilon)$ 
in such a way that the corresponding  Mori chamber 
contains the class of the adjoint divisor $[K_{\widehat{V}} + \widehat{Z}]$.   
In \cite[3.3]{Ish99}   
the simple polytope $\square(\varepsilon)$ was obtained by so called "puffing up"  
the polytope $F(P)$. We note that  Ishii's   
method does not use toric   
crepant morphisms that appeared in constructing Calabi-Yau 
minimal models corresponding to reflexive Newton polytopes $P$ and her method does not show that in general 
the minimal model $\widehat{Z}$ can be always 
obtained as nef-divisor on the simplicial toric  
variety $\widehat{V}$, i.e., that the corresponding 
Mori chamber can  be chosen  to  contain simultaneously two nef divisor classes: the adjoint class
$[K_{\widehat{V}} +  \widehat{Z}]$ and the class $[ \widehat{Z}]$.  
The latter is obtained in this paper using  the  {\em canonical 
model $\widetilde{Z}$} of $Z$ and  crepant morphisms $\widehat{Z} \to \widetilde{Z}$ that allow combinatorial computing  the Kodaira dimension $\kappa(\widehat{Z})$, the intersection number 
$(K_{\widehat{Z}})^{d-1}$ and all plurigenera of the minimal model 
$\widehat{Z}$ via the polytope $F(P)$ \cite{Gie21}.

\section{Introduction}

In this paper  we suggest a new method of "puffing up" 
the Fine interior $F(P)$ and  constructing   
minimal models of nondegenerate hypersurfaces $Z$. 
This method is a direct generalization of the construction of  
mi\-ni\-mal Calabi-Yau  
models in case of  reflexive 
Newton polytopes $P$ \cite{Bat94}. 
If $P$ is a $d$-dimensional reflexive polytope, then
its Fine interior 
$F(P)$ is just the origin $0 \in M$,  
and the support of the Fine interior 
$S_F(P)$ is exactly the set of all nonzero lattice 
points on the boundary of  the dual reflexive polytope $P^*$.

Recall that the construction of  minimal Calabi-Yau compactifications $\widehat{Z}$ of $Z$ consists of the 
following  two steps \cite{Bat94}:

\begin{itemize}
\item First,  one considers the Zariski closure $\widetilde{Z}$ 
of the affine hypersurface $Z$ in the 
Gorenstein toric Fano variety $\widetilde{V}:= V_P$ 
corresponding to the normal fan $\widetilde{\Sigma} := 
\Sigma_P$ 
of $P$. The projective hypersurface  
$\widetilde{Z} \subset \widetilde{V}$ is a Calabi-Yau variety
with at worst Gorenstein canonical singularities. From now on 
we call the 
projective hypersurface 
$\widetilde{Z}$ the {\em canonical model} of $Z \subset \T^d$. 
 \item Second,  one considers the Zariski closure $\widehat{Z}$
of $Z$ in the toric variety  $\widehat{V}$   
corresponding to a maximal projective simplicial 
refinement $\widehat{\Sigma}$ of the fan $\widetilde{\Sigma}$ satisfying the condition $\widehat{\Sigma}[1] = P^*\cap (\Z^d \setminus \{0\})$. The corresponding    
toric morphism $\widehat{V} \to \widetilde{V}$ is crepant and 
$\widehat{Z}$ is a minimal Calabi-Yau variety with at worst Gorenstein terminal singularities obtained as maximal projective
crepant partial 
desingularization (a MPCP-desingularization) of the canonical model 
$\widetilde{Z}$.  Such a desingularization is determined by choosing triangulations of facets of $P^*$, and it is not unique 
in general.  
\end{itemize}

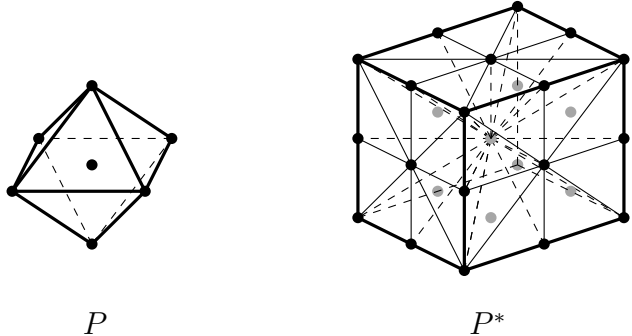
\begin{center}
\begin{tikzpicture}[scale=0.7]

\fill (0,0) circle (3pt);

\fill (1,-0.5) circle (3pt); 
\fill (-1,0.5) circle (3pt);

\fill (0,1.5) circle (3pt); 
\fill (0,-1.5) circle (3pt); 

\fill (1.5,0.5) circle (3pt); 
\fill (-1.5,-0.5) circle (3pt); 

\draw[solid,very thick] (0,1.5) --  (1.5,0.5) ; 
\draw[solid,very thick] (0,1.5) -- (1,-0.5) ;
\draw[solid,very thick] (0,1.5) -- (-1.5,-0.5) ;
\draw[solid,very thick] (0,1.5) --(-1,0.5);

\draw[solid,very thick]  (1.5,0.5) -- (1,-0.5) -- (-1.5,-0.5) -- (-1,0.5);
\draw[dashed]  (1.5,0.5) --  (-1,0.5);

\draw[solid,very thick] (0,-1.5) -- (-1.5,-0.5) ; 
\draw[solid,very thick] (0,-1.5) -- (1,-0.5) ; 

\draw[dashed] (0,-1.5) -- (1.5,0.5) ;
\draw[dashed] (0,-1.5) -- (-1,0.5) ;

\node [left] at (0.5,-3) {{\bf  $P$}}; 
\begin{scope}[xshift= 8cm] 
\fill (2,-1) circle (3pt);
\fill  (2,2) circle (3pt); 
\fill  (2,0.5) circle (3pt);
\fill  (-1,-2) circle (3pt);
\fill  (-1,-0.5) circle (3pt);
\fill  (-1,1) circle (3pt);
\fill  (-3,-1) circle (3pt);
\fill  (-3,0.5) circle (3pt);
\fill  (-3,2) circle (3pt);
\fill  (0,3) circle (3pt);
\fill  (0.5,0) circle (3pt);
\fill  (0.5,1.5) circle (3pt);
\fill  (0.5,-1.5) circle (3pt);
\fill  (-2,0) circle (3pt);
\fill  (-2,1.5) circle (3pt);
\fill  (-2,-1.5) circle (3pt);
\draw[solid] (-2,0) --  (-2,-1.5); 
\draw[solid] (-2,0) --  (-2,1.5); 

\draw[solid] (-2,0) --  (-1,-2); 
\draw[solid] (-2,0) --  (-1,-0.5); 
\draw[solid] (-2,0) --  (-1,1); 

\draw[solid] (-2,0) --  (-3,-1); 
\draw[solid] (-2,0) --  (-3,0.5); 
\draw[solid] (-2,0) --  (-3,2); 

\fill[color=gray!70!white]  (-0.5,0.5) circle (3pt); 

\draw[dashed]  (-0.5,0.5) -- (-2,-1.5); 
\draw[dashed]  (-0.5,0.5) -- (-2,1.5);
\draw[dashed]  (-0.5,0.5) -- (-1,-2); 
\draw[dashed]  (-0.5,0.5) -- (-1,-0.5); 
\draw[dashed]  (-0.5,0.5) -- (-1,1); 
\draw[dashed]  (-0.5,0.5) -- (-3,-1); 
\draw[dashed]  (-0.5,0.5) -- (-3,0.5); 
\draw[dashed]  (-0.5,0.5) -- (-3,2); 

\draw[dashed]  (-0.5,0.5) -- (2,-1); 
\draw[dashed]  (-0.5,0.5) -- (2,0.5); 
\draw[dashed]  (-0.5,0.5) -- (2,2); 
\draw[dashed]  (-0.5,0.5) -- (0.5,1.5); 
\draw[dashed]  (-0.5,0.5) -- (0.5,0); 
\draw[dashed]  (-0.5,0.5) -- (0.5,-1.5); 
\draw[dashed]  (-0.5,0.5) --  (-1,-2); 
\draw[dashed]  (-0.5,0.5) --   (-1,-0.5);
\draw[dashed]  (-0.5,0.5) -- (-1,1); 

\draw[solid] (0.5, 0) -- (2,-1); 
\draw[solid] (0.5, 0) -- (2,0.5); 
\draw[solid] (0.5, 0) -- (2,2); 
\draw[solid] (0.5, 0) --  (0.5,1.5); 
\draw[solid] (0.5, 0) -- (-1,1); 
\draw[solid] (0.5, 0) --(0.5,-1.5); 
\draw[solid] (0.5, 0) -- (-1,-2); 
\draw[solid] (0.5, 0) -- (-1,-0.5);

\fill[color=gray!70!white]  (-0.5,0.5) circle (3pt);

\fill[color=gray!70!white]  (1,-0.5) circle (3pt); 
\fill[color=gray!70!white]  (1,1) circle (3pt); 
\fill[color=gray!70!white]  (-1.5,-0.5) circle (3pt); 
\fill[color=gray!70!white]  (-0.5,-1) circle (3pt); 
\fill[color=gray!70!white]  (-1.5,1) circle (3pt); 

\fill (-1.5,2.5)  circle (3pt);  
\fill (1,2.5)  circle (3pt); 
\fill (-0.5,2)  circle (3pt);

\draw[solid] (-0.5,2) -- (-1.5,2.5); 
\draw[solid] (-0.5,2) --  (1,2.5); 
\draw[solid] (-0.5,2) -- (2,2);

\draw[solid] (-0.5,2) -- (0,3);
\draw[solid] (-0.5,2) -- (-2,1.5);
\draw[solid] (-0.5,2) -- ( -1,1);

\draw[solid] (-0.5,2) -- (-3,2);
\draw[solid] (-0.5,2) -- (0.5,1.5);

\draw[dashed] (-0.5,0.5) -- (-0.5,2); 
 \draw[dashed] (-0.5,0.5) -- (0,3); 
 \draw[dashed] (-0.5,0.5) -- (1,2.5); 
 \draw[dashed] (-0.5,0.5) -- (-1.5,2.5);

\fill[color=gray!70!white]  (0,0) circle (3pt); 
\fill[color=gray!70!white]  (0,1.5) circle (3pt); 
\fill[color=gray!70!white]  (0,0) circle (3pt);

\draw[solid, very thick] (0,3) -- (-3,2) --(-1,1) --(2,2) -- (0,3);
\draw[solid, very thick] (2,2) -- (2,-1) -- (-1,-2)  -- (-3,-1) --(-3,2) ;
\draw[solid, very thick] (-1,1) --(-1,-2);
\draw[dashed] (0,0) -- (2,-1);
\draw [dashed] (0,0) -- (-3,-1);
\draw [dashed] (0,0) -- (0,3); 

\node [left] at (0,-3) {{$ P^*$}}; 
\end{scope}
\end{tikzpicture}
\end{center}

The main idea of the new construction of minimal models $\widehat{Z}$ is based on another 
"puffing up" of the rational polytope $F(P)$ using 
the real full-dimensional polytopes
$F(\lambda P) \subset M_\R$ which are Fine interiors  
of the real $\lambda$-multiples of the lattice polytope $P$. If  $\lambda = 1 + \varepsilon$, then  
for all sufficiently small $\varepsilon > 0 $ the normal fan $\widetilde{\Sigma}$ of the full-dimensional real polytope $F((1 +\varepsilon)P)$ is independent on $\varepsilon$. Note that the fan $\widetilde{\Sigma}$ is not simplicial in general, but 
all primitive lattice vectors  $ \nu \in \widetilde{\Sigma}[1]$ belong
to the support $S_F(P)$.  Moreover, one has a Minkowski sum 
decomposition 
\[  F((1 +\varepsilon)P) = F(P) + \varepsilon C(P), \]
where 
\[  C(P) := \{ x \in M_\R \; | \; \langle x, \nu \rangle \geq \min_{p \in P} \langle p,  \nu \rangle \;\; \forall \nu \in S_F(P) \} \]
is a $d$-dimensional rational polytope containing $P$. We call 
the rational polytope $C(P)$ the {\em canonical hull of } $P$.  

The key step in our construction is to consider the toric variety 
$\widetilde{V}$ corresponding to  the normal fan $\widetilde{\Sigma}$ of the 
$d$-dimensional 
Minkowski sum 
 $$\widetilde{P} := C(P) + F(P).$$

We show that both the 
projective toric variety $\widetilde{V} = V_{\widetilde{\Sigma}}$ 
and the Zariski closure $\widetilde{Z}$ of $Z$ in $\widetilde{V}$ have at worst  $\Q$-Gorenstein canonical singularities.  Moreover,   
the canonical class $K_{\widetilde{Z}}$ is a semiample 
$\Q$-Cartier divisor on $\widetilde{Z}$. We call the projective variety 
 $\widetilde{Z}$ the 
{\em canonical model} of $Z$. 
If $\widehat{\Sigma}$ is a maximal projective simplicial 
refinement  of the fan $\widetilde{\Sigma}$ 
with the property $\widehat{\Sigma}[1] = S_F(P)$, then the corresponding toric morphism  
$\varphi\, : \; \widehat{V} \to \widetilde{V}$ is crepant, 
the toric variety $\widehat{V}$ has 
at worst terminal singularities,  
 and the Zariski 
closure $\widehat{Z}$ of $Z$ in $\widehat{V}$ is a
minimal model together with the crepant morphism 
$\varphi\, :\, \widehat{Z} \to \widetilde{Z}$.  Moreover, 
the minimal model $\widehat{Z} \subset \widehat{V}$ 
is a semiample $\Q$-Cartier divisor of $\widehat{V}$.

\medskip

The paper is organized as follows. 
\medskip

In Section 1 we deal with combinatorial tools  from convex geometry 
of $M$-lattice polytopes $P$. It turns out to be very useful to extend the class of $M$-lattice polytopes and consider 
a larger class of convex polytopes $P \subset M_\R$ which
 includes all   
rational polytopes and even some real ones. We call  polytopes 
$P$ in this class  {\em generalized Delzant polytopes} (see 
Definition \ref{Delzant}). Using the dual $N$-lattice in $N_\R$ and the natural pairing 
$\langle *,* \rangle\, :\, M \times N \to \Z$,  
we consider the associated with $P \subset M_\R$  
upper-convex 
piecewise linear function 
\[ \ord_P \;:\;  N_\R \to \R, \;\; y \mapsto \ord_P(y) := \min_{x \in P} \langle x, y \rangle \;\; (y \in N_\R),   \]
and we use it in the definitions of the following three combinatorial objects: 

\begin{itemize}
\item 
the {\em Fine interior of $P$}
\[ F(P) := \{ x \in M_\R \; | \; \langle x, \nu \rangle
\geq \ord_{P}(\nu) +1\;\; \; \forall \nu \in 
N \setminus \{0 \} \}; \]  
\item 
the {\em support of the Fine interior of} $P$ 
(if $F(P) \neq \emptyset$)
\[ S_F(P) := \{ \nu \in N \; | \; \ord_{F(P)}(
\nu)= \ord_{P} (\nu)+ 1\}; \]
\item
the   {\em canonical hull $C(P)$ of $P$} (if $F(P) \neq \emptyset$)  
\[ C(P):= \{ x \in M_\R \; | \; \langle x, \nu \rangle
\geq \ord_{P}(\nu)\;\;\;  \forall \nu \in S_F(P) \}. \]
\end{itemize}

In Section 2 we concentrate our attention  upon  the canonical hull $C(P)$ of $P$. Full-dimensional 
generalized Delzant  polytopes $P$ such that 
$P = C(P)$ we call   
{\em canonically closed}. We show that  
multiples $\lambda P$ of any full-dimensional
generalized Delzant  polytope $P$ are  always 
canonically closed for sufficiently large $\lambda \gg 1$. 

Our main interest  in Section 3  is  the relation between the canonical 
refinement $\Sigma_P^{\rm can}$ 
of the normal fan $\Sigma_P$ of  $P$ 
and its Fine interior $F(P)$. In particular, 
we show that the canonical refinement 
  $\Sigma_P^{\rm can}$ of $\Sigma_P$ is the normal 
fan of the Fine interior $F(\lambda P)$ for sufficiently 
large $\lambda \gg 1$.

Section 4 is devoted to  main combinatorial results that we use  
 in our constructing minimal models of nondegenerate hypersurfaces $Z$ via its Newton polytope $P$. We denote by 
 $\widetilde{V}$ the projective toric variety  defined by  the  fan $\widetilde{\Sigma}$ that is  the normal fan of 
  the Minkowski sum $\widetilde{P} := C(P) + F(P)$. 
The toric variety  $\widetilde{V}$ will 
play a central role in our construction.  
We show that 

1) the set $\widetilde{\Sigma}[1]$ of primitive lattice 
generators of $1$-dimensional cones of the fan 
$\widetilde{\Sigma}$ is contained in $S_F(P)$,  

2) the projective 
toric variety  $\widetilde{V}$ is $\Q$-Gorenstein and it has at worst canonical singularities, 

3) for any $2$-dimensional cone $\sigma \in \widetilde{\Sigma}(2)$ spanned by two primitive lattice vectors 
$\nu_i, \nu_j \in S_F(P)$ one has the following two properties:
\begin{itemize}
\item the $(d-2)$-dimensional affine 
linear subspace $L_\sigma^1 \subset M_\R$ defined by two linear equations 
\[ \langle x, \nu_i \rangle = \ord_P(\nu_i) +1, \;\;  
\langle x, \nu_j \rangle = \ord_P(\nu_j) +1 \]
has nonempty intersection with the Fine interior $F(P)$,

\item the $(d-2)$-dimensional affine 
linear subspace $L_\sigma^0 \subset M_\R$ defined by two linear equations 
\[ \langle x, \nu_i \rangle = \ord_P(\nu_i), \;\;  
\langle x, \nu_j \rangle = \ord_P(\nu_j) \]
contains at least one vertex of  the Newton polytope $P$.
\end{itemize}

\noindent
We note that two Minkowski summands $C(P)$ and $F(P)$ of $\widetilde{P}$ define two nef $\Q$-Cartier divisors on the toric 
variety $\widetilde{V}$, and two  
toric morphisms $\varrho\, :\, \widetilde{V} \to V_{C(P)}$
and $\vartheta\, :\, \widetilde{V} \to V_{F(P)}$ that have opposite 
properties with respect to canonical class $K$:
\[
\xymatrix{ K>0  & \widetilde{V}\ar[ld]_{\varrho} 
\ar[rd]^{\vartheta} &   K<0 \\
V_{C(P)}  \ar@{-->}[rr]
&  & V_{F(P)}}
\]

\medskip

\noindent
The toric morphism $\varrho\, :\, \widetilde{V} \to V_{C(P)}$ is birational and $K$-positive, i.e., the $\Q$-Cartier canonical 
divisor $K_{\widetilde{V}}$ has positive intersection with all $1$-dimensional 
toric stata in $\widetilde{V}$ contracted by $\varrho$.  
We show that via the toric morphism $\varrho$ one
  can consider $\widetilde{V}$ 
as minimal 
 canonical  partial resolution of $V_{C(P)}$.  The toric 
morphism 
$\vartheta\, :\, \widetilde{V} \to V_{F(P)}$ is either birational
if $\dim F(P) =d$, 
or a toric $\Q$-Fano fibration if $\dim F(P) < d$. A maximal projective  simplicial 
refinement $\widehat{\Sigma}$ of $\widetilde{\Sigma}$
with $\widehat{\Sigma}[1] = S_F(P)$ defines a crepant birational toric morphism $\varphi\, :\, \widehat{V} \to \widetilde{V}$ such that  $\widehat{V}$ 
is a projective simplicial toric variety  with 
at worst terminal singularities.

In Section 5 we define the {\em canonical
model $\widetilde{Z}$} of $Z$ as Zariski closure of $Z$ in the 
$\Q$-Gorenstein toric variety $\widetilde{V}$.  By the construction of 
$\widetilde{V}$, 
the canonical model $\widetilde{Z}$ is a big semiample $\Q$-Cartier divisor on $\widetilde{V}$ corresponding to the polytope $C(P)$. 
The  crucial is the fact that 
the canonical model $\widetilde{Z}$ is a normal variety   
satisfying 
the adjunction formula, so that the canonical class of $\widetilde{Z}$ is the restriction of the adjoint nef $\Q$-Cartier 
divisor $K_{\widetilde{V}} + \widetilde{Z}$ on the toric variety $\widetilde{V}$ to the semiample hypersurface $\widetilde{Z}$.

In Section 6  we show that  the 
hypersurface $\widetilde{Z} \subset \widetilde{V}$ 
has at worst $\Q$-Gorenstein canonical singularities and  
we obtain minimal models $\widehat{Z}$ 
of $Z$ as Zariski closures of $Z \subset \T^d$ in the simplicial terminal toric varieties $\widehat{V}$ defined by the simplicial 
refinements  $\widehat{\Sigma}$ of the normal fan $\widetilde{\Sigma}$ satisfying 
the condition $\widehat{\Sigma}[1] = S_F(P)$. 

In Section 7,  using the toric morphism $\vartheta\, :\, \widetilde{V} \to V_{F(P)}$, we investigate the Iitaka fibration of $\widetilde{Z}$
defined by the pluricano\-nical ring  
\[ R:= \bigoplus_{m \geq 0} H^0(\widetilde{Z}, 
m K_{\widetilde{Z}}).  \] 
In particular, 
we prove the following formula for the Kodaira dimension  
\[ \kappa(\widetilde{Z}) = \min \{ \dim F(P), d-1 \}, \]
and compute the intersection number $(K_{\widetilde{Z}})^{d-1}$:
$$(K_{\widetilde{Z}})^{d-1} = \begin{cases} 
{\rm Vol}_d(F(P)) + \sum_{ Q \prec F(P) \atop \dim Q =d-1} {\rm Vol}_{d-1}(Q), & \text{ if $\dim F(P) =d$;} \\
2 {\rm Vol}_{d-1}(F(P)), & \text{ if $\dim F(P) =d-1$;} \\
0  , & \text{ if $\dim F(P) < d-1$,}
\end{cases}
$$
where ${\rm Vol}_k(\cdot)=k!vol(\cdot)$ denotes the $k$-dimensional volume normalized by the condition ${\rm Vol}_k(\Delta_k) =1$ if 
$\Delta_k$ is the basic $k$-dimensional lattice simplex. 

Recently  Julius Giesler found  
a nice general combinatorial formula  for 
all plurigenera $g_m(\widetilde{Z})$, $m \geq 0$  \cite{Gie21}. 
His formula yields another 
method of proving the above combinatorial formula for 
$(K_{\widetilde{Z}})^{d-1}$. 

It is important fact that 
the intersections of the canonical hypersurface 
$\widetilde{Z} \subset \widetilde{V}$ with the $(d- \dim F(P))$-dimensional 
generic fibers of  toric morphis $\vartheta\, : \, \widetilde{V} \to 
V_{F(P)}$ are nondegenerate toric hypersurfaces of Kodaira dimension $0$ embedded into 
canonical toric $\Q$-Fano 
varieties of dimension $d - \dim F(P)$.  In particular, 
we obtain that the toric pairs
$(\widetilde{V}, \widetilde{Z})$ together with the toric 
morphism 
$\vartheta\, :\, \widetilde{V} \to V_{F(P)}$ corresponding to   
the adjoint nef $\Q$-Cartier 
divisor $K_{\widetilde{V}} + \widetilde{Z}$ 
provide toric  examples of the log Calabi-Yau fibrations suggested recently 
by Caucher Birkar \cite{Bir18b,Bir21}. 
\medskip

The last section 8 gives a short overview of some possible further developments of the considered topic and contains open questions. 
\medskip

{\bf Acknowledgements.} I would like to express my thanks for 
helpful discussions to Martin Bohnert, Dimitios Dais, Jonatan Fine, Julius Giesler, J\"urgen Hausen, Shihoko Ishii, Alexander Kasprzyk, Ivo Radloff, Miles Reid, Karin Schaller, Harald Skarke, and especially to Jaron  Treutlein for his constant interest towards interactions between MMP and lattice polytopes \cite{Tre10}.

\newpage
 
\section{Main combinatorial objects}

Let $M \cong \Z^d$ be a lattice of rank $d$, $N:= \Hom(M, \Z)$
the dual lattice. Denote by 
 $\langle *, * \rangle \,:\,  M \times N \to \Z$ 
the natural pairing. We extend this pairing to the real 
vector spaces 
$M_\R$ and $N_\R$. 
For any convex compact subset 
$\square \subset M_\R$, we consider the upper-convex 
real valued 
function $\ord_\square 
\, :\, 
N_\R \to \R$ defined as $\ord_\square(y) :=  \min_{x \in \square} \langle x, y \rangle $.

The main combinatorial object in our study is a full-dimensional convex 
lattice polytope  $P \subset M_\R$,  i.e., 
a $d$-dimensional 
convex hull of a finite set $A \subset M$.  We denote by 
$\Sigma_P$ the normal fan of $P$.  
  There is a bijection between $k$-dimensional 
  faces $Q \preceq P$ $(0 \leq k \leq d)$ and $(d-k)$-dimensional 
cones 
$$\sigma_Q := \{ y \in N_\R \mid  \ord_P(y) = 
\< x, y \>   \; \forall x \in Q \}  \in \Sigma_P$$  
in the dual space $N_\R$. The convex function $\ord_P \, : \, N_\R \to \R$ is $\Sigma_P$-piecewise linear,  it  defines an ample Cartier divisor ${\mathcal L}_P$ on the corresponding toric variety 
$V_P := V_{\Sigma_P}$. Obviously,  we have 
\[  \ord_P(y) = \min_{m \in P\cap M} \langle m, y \rangle  = \min_{m \in A}
\langle m, y \rangle. \]

It turns out that our study of $d$-dimensional lattice polytopes requires  
a larger class of convex polytopes
under consideration. First of all, we need to  work with 
rational polytopes in $M_\R$ 
having dimension $\leq  d$. 
Moreover, we must allow to multiply any   
considered convex 
polytope by positive real numbers and to shift these polytopes by real 
vectors  $x \in M_\R$. At the same time, we want to keep 
connections  of considered 
convex polytopes with the lattices $M$ and $N$.  
All these requirements lead us to  a larger  class of convex polytopes  of  
dimension $\leq d$ that we call {\em generalized Delzant
polytopes}. This notation  is inspired by the notion of Delzant  polytopes
in the symplectic toric geometry \cite{Del88,Gui94,daS01}:

\begin{defi} \label{Delzant}
Let $P \subset M_\R$ be a compact convex polytope 
of dimension $\leq d$. We call $P$ a {\em generalized Delzant polytope} if 
there exists a finite set of lattice vectors  
${\nu_1, \ldots, \nu_s} \subset N$ and 
a finite set of real numbers $\delta_1, \ldots, \delta_s \in \R$ such that 
\[ P = \{ x \in M_\R \mid  \langle x, \nu_i \rangle \geq \delta_i, \; \; 
1 \leq i \leq s  \}. \] 
The normal fan  $\Sigma_P$ of a generalized 
Delzant polytope $P$
is a 
generalized fan in $N_\R$  defining a toric variety $V(\Sigma_P)$ \cite[Def. 6.2.2]{CLS11} and 
the $\Sigma_P$-piecewise linear function $\ord_P$
defines an ample $\R$-Cartier divisor ${\mathcal L}_P$ on 
$V(\Sigma_P)$. 
\end{defi} 

\begin{ex}
A two-dimensional generalized Delzant polytope
 $P \subset \R^2$ is simply 
a convex $n$-gon whose sides have rational slopes. 
\end{ex}

\begin{defi} 
Let $P \subset M_\R$ be a full dimensional 
generalized Delzant polytope. We 
define the {\em Fine interior $F(P)$} of $P$ 
as 
\[  F(P):= \{ x \in M_\R \mid \langle x, \nu \rangle \geq \ord_P(\nu) +1 
\;\; \forall \nu \in N \setminus \{0 \} \}. \] 
\end{defi}

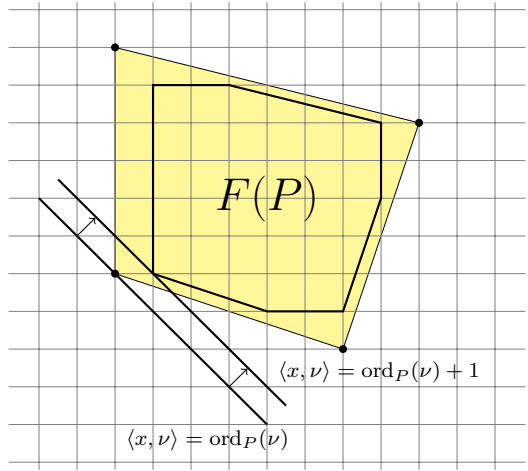
\begin{center}
\begin{tikzpicture}
\draw[ fill=gray!20!white, very thin] (-2,-1) -- (-2,2) -- (2,1) -- (1,-2)
-- ( -2,-1) ;
\draw[ fill=gray!45!white,thick]  (-1.5,-1) -- (-1.5,1.5) -- (-0.5,1.5) -- (1.5,1) --
(1.5, 0) -- (1, -1.5) -- (0,-1.5) -- (-1.5, -1);
\draw[step=0.5cm,gray,very thin]
(-3.4,-3.6) grid (3.4,2.6) ;

\fill (-2,2) circle (1.5pt);
\fill (2,1) circle (1.5pt);
\fill (1,-2) circle (1.5pt);
\fill (-2,-1) circle (1.5pt);
\draw[thick] (-3,0) -- (0,-3); 
\draw[thick] (-2.75,0.25) -- (0.25,-2.75); 
\draw[ - >] (-2.5,-0.5) -- (-2.26,-0.26);
\draw[ - >] (-0.5,-2.5) -- (-.26,-2.26);
\node at (0,0) {\Large {\bf $F(P)$}};
\node[right]  at (-2,-3.2)  {\tiny  {\bf $\langle x, \nu \rangle = {\ord}_P(\nu)$}};
\node[right]  at (0,-2.3)  {\tiny  {\bf $\langle x, \nu \rangle = {\ord}_P(\nu) + 1$}}; 

\end{tikzpicture}
\end{center}

\begin{rem} \label{fp-shift}
Obviously, the Fine interior $F(P)$ is a compact 
subset strictly contained  in the interior $P^\circ$ of $P$ and $F(P)$ commutes with shifts by any real vectors,  i.e., for any generalized Delzant polytope $P$ one has
\[ F( P + x) = F(P) + x, \;\; \forall x \in M_\R. \]
\end{rem}

\begin{rem} \label{int-latt}
Let  $P$ be a full dimensional lattice polytope. Since $\ord_P(N) \subset \Z$, every interior lattice point 
$m \in P^\circ \cap M$ is 
contained in $F(P)$.  So we obtain the inclusion 
 ${\rm Conv}(P^\circ \cap M) \subseteq F(P)$. In particular,  
 $F(P) \neq \emptyset$ as soon as  $P$ contains at least one interior lattice point.  If $\dim P = 2$, one has the equality 
\[ F(P) = {\rm Conv}(P^\circ \cap M), \]
see \cite[Prop. 2.9]{Bat17}. This equality does not hold in general 
if $\dim P \geq 3$. 
\end{rem}

\begin{ex}  \label{exam-parall}
Let $P$ be a $3$-dimensional parallelepiped with $8$ lattice vertices 
\[ (0,0,0), \, (-1,1,1),\, (1,-1,1), \, (1,1,-1),\, (2,0,0), \,(0,2,0),\, (0,0,2),\, 
(1,1,1). \]
Then $P$ has no interior lattice points, but  $F(P) = 
\{ (1/2,1/2,1/2 ) \} \neq \emptyset$. Note that  the Fine interior of the rational shifted 
polytope $P'= P -  (1/2,1/2,1/2 )$ is $\{ 0 \}$, and the polar dual 
of $P'$ is the lattice polytope $${\rm Conv}\{ \pm (0,1,1), \pm (1,0,1), \pm (1,1,0) \}.$$
\end{ex}

\begin{rem} \label{fp-monoton}
Note that the Fine interior is {\em monotone}, i.e., 
\[ P \subseteq P' \Rightarrow F(P ) \subseteq F(P'), \]
because $\ord_P(\nu) \geq \ord_{P'}(\nu)$ for all $\nu \in N$. 
In particular, one has $F(P) \subseteq 
F(\lambda P)$ for any real $\lambda \geq 1$ if $P$ 
contains the origin $0 \in M$.  Using shifts \ref{fp-shift}, one easily obtains  
\[ F(P) \neq \emptyset \Rightarrow F(\lambda P) \neq \emptyset 
\;\; \forall \lambda  \geq 1. \] 
\end{rem}

\begin{prop}
Let $P$ be a full dimensional generalized Delzant polytope 
$P \subset M_\R$. There exists a positive real number $\lambda_P^0$ such that $F(\lambda P) \neq \emptyset$ if and 
only if $\lambda \geq \lambda_P^0$. Moreover, the polytope 
$F(\lambda P)$ is full dimensional if and only if $ \lambda > \lambda_P^0$. 
\end{prop}

\noindent
{\bf Proof.} By \ref{fp-shift}, we may assume that $P$ contains the 
origin $0 \in M$. We set 
\[ \lambda_P^0 := \inf \{ \lambda \in \R_{>0} \mid F(\lambda P) 
\neq \emptyset \}. \]
By \ref{fp-monoton}, $F(\lambda_P^0 P) \neq \emptyset$, because
$M_\R$ is a complete Banach space.  For all  
$\lambda > \lambda_P^0$ the polytope $F(\lambda P)$ contains the full dimensional Minkowski sum 
$F(\lambda_P^0 P) + (\lambda - \lambda_P^0) P$.  Thus,  $F(\lambda P)$ is full 
dimensional for all  
$\lambda > \lambda_P^0$. On the other hand, $F(\lambda_P^0 P)$ can not be full dimensional, because otherwise $F((\lambda_P^0 - \varepsilon)P)$ would  be full dimensional for sufficiently small $\varepsilon >0$ and the latter contradicts 
$F((\lambda_P^0 - \varepsilon)P) = \emptyset$. 
\hfill  $\Box$

\begin{defi} \label{supportFP}
Let $P$ be a $d$-dimensional generalized Delzant polytope. Assume that  $F(P) \neq \emptyset$.  
A nonzero lattice point $\nu \in N$ is called {\em essential} for 
$F(P)$ if 
\[ \ord_{F(P)}(\nu) = \ord_P(\nu) + 1. \]
The set of all essential lattice points $\nu \in N$ for $F(P)$ 
is called  the {\em support of $F(P)$} and we denote 
it by $S_F(P)$.  
\end{defi}

\begin{rem} \label{prop-supportFP}
Using  \ref{fp-shift} , we obtain  
that $S_F(P+x) = S_F(P)$ for all $x \in M_\R$. 
We call  a lattice vector $\nu \in N$ in a cone $\sigma \in \Sigma_P$ of the normal fan 
{\em  $\Sigma_P$-irreducible} if $\nu \neq  \nu_1 + \nu_2$ for some $\nu_1,\nu_2 \in \sigma \cap N$. Since the function $\ord_P$ is linear on every cone $\sigma \in \Sigma_P$, we obtain that two inequalities  
 \[ \langle x, \nu_1 \rangle \geq \ord_P(\nu_1) +1,  \;\; \langle x, \nu_2 \rangle \geq \ord_P(\nu_2) +1 \]
for $\nu_1, \nu_2 \in \sigma \cap N$ 
imply $\langle x, \nu \rangle \geq \ord_P(\nu) +2 > \ord_P(\nu) +1$, where $\nu = \nu_1 + \nu_2$, i.e. $\nu \not\in S_F(P)$.  Therefore, every $\nu \in S_F(P)$ must be $\Sigma_P$-irreducible. By Gordan's lemma, the set of $\Sigma_P$-irreducible elements $
\nu \in N$ for any $\sigma  \in \Sigma_P$ is finite (Hilbert basis). So  $S_F(P)$ is a finite set consisting of some $\Sigma_P$-irreducible primitive lattice vectors $\nu \in N$  
and the Fine interior $F(P)$ is a generalized Delzant polytope of dimension $\leq d$. Moreover, one has
\[ F(P) =  \{ x \in M_\R \mid \langle x, \nu \rangle \geq \ord_P(\nu) +1 
\;\; \forall \nu \in S_F(P) \}, \]
because the compact polytope $F(P)$ is strictly contained in the open 
halfspace   $$\{ x \in M_\R \mid\langle x, \nu \rangle > \ord_P(\nu) +1\}, \;\; {\rm if} \; \nu \not\in S_F(P),$$ and we can ignore 
the corresponding inequiality $\langle x, \nu \rangle \geq \ord_P(\nu) +1 $ for $F(P)$. 
 \end{rem}
\medskip

Note that computer calculations of the Fine interior $F(P)$ and its support $S_F(P)$ that are using the Hilbert basis of lattice semigroups $\sigma \cap N$ can be rather time 
consuming. Another way for such a calculation of $F(P)$ and $S_F(P)$ is provided by  the next useful statement: 

\begin{prop} \label{Fp-computer}
Let $P$ be a full dimensional generalized Delzant polytope. 
Denote by $\Sigma_P[1] \subset N$ 
the set of all primitive   
inward pointing facet normals of $P$. Then  
\[ S_F(P) \subseteq {\rm Conv}(\Sigma_P[1]). \] 
In particular,  the set $S_F(P)$ is finite  and 
\[ F(P) = \{ x \in M_\R \mid \langle x, \nu \rangle \geq \ord_P(\nu) +1 
\;\; \forall  \nu \in {\rm Conv}(\Sigma_P[1]) \cap 
\{ N \setminus \{0\} \}. \] 
\end{prop}

\noindent
{\bf Proof.} First, we note that $0 \in {\rm Conv}(\Sigma_P[1])$, because the lattice vectors in  
 $\Sigma_P[1]$ generate the complete normal  
 fan $\Sigma_P$. 
Take a lattice vector $\nu  \in S_F(P)$. Since the normal 
fan  $\Sigma_P$ is complete, $\nu$  is contained in some $d$-dimensional cone $\sigma \in \Sigma_P(d)$.  
Consider the set of primitive lattice vectors     
$$\{ \nu_1, \nu_2, \ldots, \nu_k\}  := 
\sigma \cap \Sigma_P[1] $$
such that $\sigma = \sum_{i=1}^k \R_{\geq 0} \nu_i$.  
It is sufficient 
to show that $\nu \in {\rm Conv}\{0, \nu_1, \ldots, \nu_k \}$, 
because ${\rm Conv}(0, \nu_1, \ldots, \nu_k) \subset {\rm Conv}(\Sigma_P[1])$. 
 
Assume that $\nu \not\in {\rm Conv}(0, \nu_1, \ldots, \nu_k)$. It follows from the inclusion $\R_{\geq 0} \nu \subset \sigma$ that  
there exists a positive $\lambda \in \R$ such that 
$\lambda < 1$ and $\lambda \nu \in  
{\rm Conv}(\nu_1, \ldots, \nu_k)$. In particular, one can write 
$\lambda \nu$ as a nonnegative  linear combination $\sum_{i=1}^k \lambda_i \nu_i$ 
with 
$\sum_{i=1}^k \lambda_i = 1$.  
Take any $x \in F(P)$ and consider 
the  linear combination of $k$ inequalities 
$ \langle x, \nu_i \rangle \geq \ord_P(\nu_i) +1$
with the coefficients $\lambda_i$ $(1 \leq i \leq k)$.  We obtain 
$ \langle x, \lambda \nu \rangle  \geq \ord_P(\lambda \nu) + 1$, 
because $\ord_P$ is linear on $\sigma$.  So we get   
 \[ \langle x, \nu \rangle \geq \ord_P(\nu) + \lambda^{-1} > 
 \ord_P(\nu) + 1 \;\;\; \forall x \in F(P). \]
 Hence $\nu \not\in S_F(P)$. Contradiction. 
 \hfill  $\Box$

\begin{ex} 
Let  $P$ be the $3$-dimensional lattice parallelepiped from 
\ref{exam-parall}. Then $S_F(P)$ consists of $6$ lattice 
vectors 
\[ S_F(P) = \{ \pm (0,1,1), \pm (1,0,1), \pm (1,1,0) \}. \]
\end{ex}

\begin{defi}
Let $P$ be a $d$-dimensional generalized Delzant polytope 
with $F(P) \neq \emptyset$. We call the polytope 
\[ C(P):= \{ x \in M_\R \mid 
 \langle x, \nu \rangle \geq \ord_P(\nu) \;\; 
\forall \nu \in 
S_F(P) \} \]
{\em the canonical hull of $P$}.  
Note that $C(P)$ contains $P$, because the inequality 
$\langle x, \nu 
\rangle \geq \ord_P(\nu)$ holds  for all $x \in P$ and all $\nu \in N$. 
In particular, $C(P)$ is a generalized Delzant polytope of dimension $d$. The compactness of $C(P)$  follows from the compactness of $F(P)$ and from the equality 
\[ F(P) =  \{ x \in M_\R \mid \langle x, \nu \rangle \geq \ord_P(\nu) +1 
\;\; \forall \nu \in S_F(P) \} \]
in Remark \ref{prop-supportFP} showing that the cone 
spanned by $S_F(P)$ equals  $N_\R$. 
\end{defi}
\medskip

\begin{rem}
If $P$ is a $d$-dimensional lattice  polytope with $F(P) \neq \emptyset$, then $C(P)$ is 
a rational polytope. In general, $C(P)$  is 
not a lattice polytope. 
\end{rem}

\begin{ex} \label{can-rat}
Let $P \subset \R^5$ the $5$-dimensional lattice polytope 
defined by the inequalities 
$$x_i \geq 0 \; (1 \leq i \leq 5), \;\; 
 x_1 + x_2 + x_3 + x_4 + 2 x_5 \leq 7, \;\;x_5 \leq 3. $$
Then $F(P) = \{(1,1,1,1,1) \}$ and 
\[ S_F(P) = \{ e_1,e_2,e_3, e_4, e_5, 
-e_1-e_2 -e_3 -e_4 -2e_5   \}. \]  
The canonical hull $C(P)$ is a rational $5$-simplex containing  the rational vertex $(0,0,0,0,7/2)$.  
Note that the primitive lattice inward normal vector to the facet 
$\{ x_5  =3 \}$ of $P$ is not an element of $S_F(P)$.  
\end{ex}

\begin{prop} \label{product}
Let $P \subset M_\R$ and $P' \subset M_R'$ be two 
full dimensional generalized Delzant polytopes 
with $F(P) \neq \emptyset$ and $F(P') \neq \emptyset$. Then 
the Cartesian 
product $$P \times P' \subset M_\R \times M_\R'$$ is 
a full dimensional generalized Delzant polytope 
with $F(P \times P') \neq \emptyset$ 
and the following statements hold: 

{\rm (a)} $F(P\times P') = F(P) \times F(P')$; 

{\rm (b)} $S_F(P\times P') = \{ (\nu, 0)  \mid
\nu \in S_F(P) \}  \cup \{ (0, \nu')  \mid \nu' \in S_F(P') \}$; 

{\rm (c)} $C(P\times P') = C(P) \times C(P')$. 

\end{prop}

\noindent
{\bf Proof.} If 
$$P= \{ x \in M_\R \mid \langle x, \nu_i \rangle \geq 
\delta_i, \;\; (1 \leq i \leq k)\} , $$ 
$$P'= \{ x \in M_\R \mid \langle x, \nu_j' \rangle \geq 
\delta', \;\; (1 \leq j \leq l) \}, $$ then the product $P \times P'$ is defined by the 
inequalities 
 \[ P \times P' = \{ (x, x') \in M_\R \times M_\R' 
\mid \langle (x,x'), (\nu_i, 0) \rangle \geq 
\delta_i\;   \forall i,  \;  \langle (x,x'), (0, \nu_j') \rangle \geq 
\delta_j', \;  \forall  j. \}\]
So $P \times P'$ is also a generalized Delzant polytope. 
For any $(\nu, \nu')\in N \times N'$ we have 
\[ \ord_{P \times P'} (\nu, \nu') = \ord_P (\nu) + \ord_{P'}(\nu'). \]
(a) If $(x, x') \in F(P \times P')$,  then 
\[  \langle (x,x') , (\nu_i, 0) \rangle \geq \ord_{P \times P'} (\nu, 0) 
\geq \ord_{P\times P'}((\nu, 0)) + 1 \;\;  
\forall \nu \in N \setminus \{0 \},  \]  
\[  \langle (x,x'), (\nu_i, 0) \rangle \geq \ord_{P \times P'} (0, \nu') 
\geq \ord_{P\times P'}((0, \nu')) + 1 \;\;  
\forall \nu' \in N' \setminus \{0 \},  \]
i.e., $(x, x') \in F(P) \times F(P')$.  
In follows that $F(P \times P') \subseteq F(P) \times F(P')$.
On the other hand, if  $(x, x') \in  F(P) \times F(P')                             $ and $(\nu, \nu') \neq (0,0)$, then 
\[    \langle (x,x'), (\nu, \nu') \rangle =  \< x, \nu \>    + 
 \< x',  \nu' \> \geq  \ord_P(\nu) + \ord_{P'}(\nu') + 1 \geq 
 \ord_{P \times P'}( \nu, \nu') +1. \]  
(b) The equality  
\[ \ord_{P \times P'} (\nu, \nu') +1 = \ord_{F(P) \times F(P')}
(\nu, \nu') \]
can hold only if $\nu =0$ or $\nu' =0$. The latter implies 
$\nu' \in S_F(P')$ or $\nu \in S_F(P)$. \\
(c) follows from (b).
\hfill $\Box$

\begin{prop} \label{can-ord}
Let $C(P)$ be the canonical hull of a generalized Delzant polytope $P$ with $F(P) \neq \emptyset$. Then 

{\rm (a)} $\ord_P(\nu) = \ord_{C(P)}(\nu)$ for all 
$\nu \in S_F(P)$; 

{\rm (b)}  $F(P) = F(C(P))$, i.e., the polytopes $P$ 
and $C(P)$ have the same 
Fine interior. 
\end{prop}

\noindent
{\bf Proof. } 
(a) It follows from 
\[ C(P) : = 
\{ x \in M_\R \mid  \langle x, \nu \rangle \geq \ord_P(\nu) \;\; 
\forall \nu \in S_F(P) \} \]
that $\ord_{C(P)}(\nu) \geq \ord_P(\nu)$  for all 
$\nu \in S_F(P)$.  The opposite inequality $\ord_P(\nu) \geq  \ord_{C(P)}(\nu)$ holds for all $\nu \in N$,  
because $P \subseteq C(P)$. 

(b) By Remark \ref{prop-supportFP},   
\[ F(P) = 
\{ x \in M_\R \mid   \langle x, \nu \rangle \geq \ord_P(\nu) +  1 \;\; 
\forall \nu \in S_F(P) \}. \]
On the other hand, we have  
the inclusion  
\[  F(C(P)) \subseteq 
\{ x \in M_\R \mid   \langle x, \nu \rangle \geq \ord_{C(P)}(\nu) +  1 \;\; 
\forall \nu \in S_F(P) \}.   \]
By  (a), we have $\ord_P(\nu) = \ord_{C(P)}(\nu)$ for all 
$\nu \in S_F(P)$. This implies the inclusion $ F(C(P))  \subseteq F(P)$.

Since $P$ is contained in $C(P)$ we obtain the opposite inclusion 
$F(P) \subseteq F(C(P))$. 

\hfill $\Box$
  
\begin{coro} \label{can-suppF}
The support of the Fine interior of $C(P)$ 
equals $S_F(P)$. 
\end{coro}

\noindent
{\bf Proof. }  All lattice vectors  $\nu \in S_F(P)$ are contained in the support of the Fine interior of $C(P)$, because 
\[  \ord_{F(P)}(\nu)  - \ord_{C(P)}(\nu) =  \ord_{F(P)}(\nu) - 
\ord_P(\nu) =1  
 \; \; \; \forall \nu \in S_F(P).   \]
The inclusion $P \subseteq C(P)$ implies that 
\[ \ord_{F(P)}(\nu)  - \ord_{C(P)}(\nu)  \geq \ord_{F(P)}(\nu) - 
\ord_P(\nu) \;\;\; \forall \nu \in N. \]
Therefore, it follows from $\ord_{F(P)}(\nu)  - \ord_{C(P)}(\nu) =1$
that $\ord_{F(P)}(\nu) - 
\ord_P(\nu) \leq 1$. On the other hand,  $\ord_{F(P)}(\nu) - 
\ord_P(\nu) \geq 1$ holds for all $\nu \in N \setminus \{0\}$. 
So, we obtain 
 $\ord_{F(P)}(\nu) - 
\ord_P(\nu) = 1$, and the support of  the Fine interior of $C(P)$ is contained
in $S_F(P)$. 
\hfill $\Box$

\begin{coro} \label{CC=C}
Let $P$ be a full dimensional generalized Delzant polytope with 
$F(P) \neq \emptyset$. Then 
\[ C(C(P)) = C(P). \]
\end{coro}

\noindent
{\bf Proof. } The statement follows straightforward  
from \ref{can-suppF} and \ref{can-ord}.  \hfill $\Box$ 
\medskip

\section{Canonically closed polytopes}

\begin{defi}
Let $P$ be a  full dimensional generalized 
Delzant  polytope $P \subset M_\R$ 
with $F(P) \neq \emptyset$. We call 
$P$ {\em canonically closed} if $C(P) = P.$
\end{defi}

\begin{rem}
It follows from \ref{CC=C} that the canonical hull $C(P)$ of a 
full dimensional  
Delzant polytope $P$ is always canonically closed. 
\end{rem}

\begin{prop} \label{can-closed-crit}
Let $P$ be a full dimensional generalized Delzant polytope 
with  $F(P) \neq \emptyset$.  Then 
$P$ is canonically closed if and only if for any facet $Q \prec P$ the primitive inward pointing    
facet normal $\nu_Q$ belongs to $S_F(P)$, i.e., its satisfies the equation $\ord_{F(P)}(\nu_Q) = \ord_P(\nu_Q) + 1$.  
\end{prop}

\noindent
{\bf Proof.} Assume that  $P$ is canonically closed. Then 
\[ P = C(P) =  \{ x \in M_\R \mid  \langle x , v 
\rangle \geq \ord_P(\nu),  \;\;\; \forall \nu \in 
S_F(P) \}. \]
Therefore, any facet $Q$ of $P$ is defined by the  
equation $\langle x , \nu 
\rangle = \ord_P(\nu)$ for some $\nu \in S_F(P)$. The primitive inward pointing  lattice normal vector $\nu_Q$ 
to a given facet $Q$ 
is uniquely determined by $Q$.  So  $\nu_Q \in S_F(P)$. 

On the other hand, one can always write  a convex 
full dimensional polytope $P$ as intersection 
of all half-spaces $\langle x , \nu_Q 
\rangle \geq  \ord_P(\nu_Q)$ corresponding to its facet $Q \prec P$: 
\[ P = \bigcap_{Q \prec P \atop \dim Q = d-1}  \{ x \in M_\R \mid \langle x , \nu_Q 
\rangle \geq \ord_P(\nu_Q) \}. \]
If all $\nu_Q$ belong to $S_F(P)$, then $C(P) \subseteq P$. 
The opposite inclusion $P \subseteq C(P)$ holds for all $P$.
Thus, we obtain $P = C(P)$. 
\hfill $\Box$  

\begin{prop} \label{cp-d=2}
Let $P$ be a $2$-dimensional  lattice polytope with
$F(P) \neq \emptyset$. Then 
$C(P) = P$.
\end{prop}

\noindent
{\bf Proof.}  By \ref{can-closed-crit}, it suffices to show that the 
primitive 
inward pointing  lattice  normal   
$\nu_Q$ to any $1$-dimensional 
side $Q \prec P$ belongs to $S_F(P)$. 
Take a $1$-dimensional 
side $Q \prec P$. By \ref{int-latt}, we have  
$F(P) = {\rm Conv}(P^\circ \cap M)$. In particular, 
$\{P^\circ \cap M\} \neq \emptyset$ and we can find a lattice 
point $m \in  F(P)$ such that $\ord_{F(P)}(\nu_Q) = 
\langle m, \nu_Q \rangle$. Take two lattice points $m_1, m_2 \in Q$
such that $m_1 - m_2$ is a primitive lattice vector. 
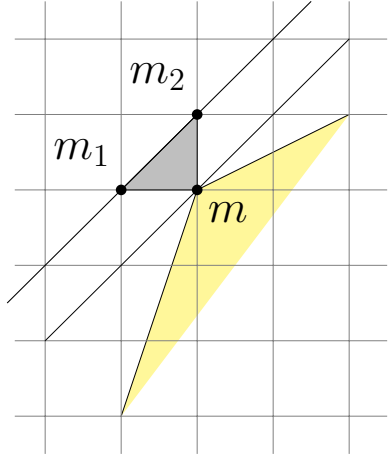
\begin{center}
\begin{tikzpicture}[scale = 1.0]


\draw[ fill=gray!30!white, very thin] (-1,-3) -- (0,0) -- (2,1) ;

\draw[step=1cm,gray,very thin]
(-2.4,-3.5) grid (2.5,2.5) ;

\draw[ fill=gray!50!white, very thin] (-1,0) -- (0,0) -- (0,1) --
(-1,0);

\fill (0,0) circle (2pt);
\fill (-1,0) circle (2pt); 
\fill (0,1) circle (2pt); 

\draw[ fill=gray!50!white, very thin] (-2.0,-2) -- (2.0, 2.0); 
\draw (-2.5,-1.5) -- (1.5, 2.5); 
\node[left] at (0, 1.2) {{ $m_2$}};
\node[left] at (-1, 0.2) {{$m_1$}};
\node[right] at (0, -0.3) {{$m$}};
\end{tikzpicture}
\end{center}

Then 
the lattice triangle $\tau:= 
{\rm Conv}(m, m_1, m_2)$ is a standard basic 
 triangle, because it follows from 
\[    \langle m, \nu_Q \rangle = \min_{a \in F(P) \cap M} \langle a, \nu_Q \rangle \]
that all lattice points in $\tau$ are only its vertices 
$m, m_1, m_2$. 
Therefore, the lattice distance between the segment 
$[m_1,m_2]$ and $m$ is $1$, and we obtain   
\[  \langle m, \nu_Q \rangle = 
\langle m_1, \nu_Q \rangle +1= \langle m_2, \nu_Q \rangle +1. \]
Thus, $\nu_Q \in S_F(P)$. 
\hfill $\Box$

\begin{rem}
The statement in \ref{cp-d=2} does not hold true 
in general for lattice polytopes $P$ of dimension $d \geq 3$, 
or for $2$-dimensional rational polytopes. 
\end{rem}

\begin{ex}
Let   $P := {\rm Conv}\{ e_0, e_1,e_2, e_3 \} \subset M_\R$ be a 
$3$-dimensional lattice simplex with vertices  
 \[ e_0 := (-1,-1,-1), \; e_1 := (1,1,0), \; e_2 := (1,0,1), \; 
 e_3 := (0,1,1). \]
 Then $F(P) = \{0 \}$ and $C(P) ={\rm Conv}\{ e_0, e_1,e_2, e_3, 
 e_4 \}$, where $e_4 := (1,1,1)$,   
i.e., $C(P) \neq P$. Note that the  inward pointing 
facet normal $\nu_Q = (-1,-1,-1) \in N$ to the facet 
$Q := {\rm Conv}\{ e_1,e_2,e_3 \} 
\prec P$ is not an element of $S_F(P)$, because $Q$ has lattice 
distance $2$ from the origin $0 = F(P)$.  
\end{ex}

\begin{ex} Let $P \subset \R^2$ be the triangle with 
vertices $(-1,0), (0,3/2), (4, -5/2)$. Then the Fine interior  
$F(P)$ is the triangle with vertices $(0,0), (0,1/2), 
(1,-1/2)$. The canonical hull $C(P)$ is strictly larger
than the triangle $P \subsetneq C(P) $, i.e., 
$P$ is not canonically closed. The integral 
distance between the side $Q:= [ (-1,0), (0,3/2)]$ and 
the rational vertex $(0,1/2) \in F(P)$ is 
$2$. The canonical hull $C(P)$ is a quadrilateral  containing 
 one more vertex $(-1,1/2)$.  
\end{ex}

\begin{center}
\begin{tikzpicture}[scale = 1.0]


\draw[ fill=gray!20!white, very thin] (-1,0) -- (0,3/2) -- (4,-5/2) 
-- ( -1,0) ;

\draw[ fill = gray!35!white, very thin] (-1,0) -- (-1, 1/2) -- (0,3/2);

\draw[fill = gray!55!white, very thin] (0,0) -- (0, 0.5) -- (1, -0.5) -- (0,0) ; 

\draw[step=1cm,gray,very thin]
(-2.4,-3.5) grid (4.5,2.5) ;

\fill (0,0) circle (2pt);
\fill (0,3/2) circle (2pt);
\fill (-1,0) circle (2pt); 
\fill (4,-5/2) circle (2pt); 
\fill (-1,1/2) circle (2pt); 

\fill (0,0.5) circle (2pt);
\fill (1,-0.5) circle (2pt); 


\node[left]  at (-1,1/2)  {{(-1,1/2)}};
\node[right] at (4, -5/2) {(4,-5/2)};
\node[right] at (0, 3/2) {(0,3/2)};
\node[left] at (-1, -1/2) {(-1,0)};
\node[left] at (0.2, -0.3) {(0,0)};
\end{tikzpicture}
\end{center}

\begin{defi} \cite{Bat94} Let $P \subset M_\R$ be 
a $d$-dimensional lattice polytope  containing 
the origin  $0 \in M$ in its interior. The polytope 
$P$  
is called {\em reflexive} if the dual 
polytope 
\[ P^* := \{ y \in N_\R \mid  \langle x, y \rangle \geq -1, 
\;\; \forall x \in P \} \] 
is a lattice polytope. 
\end{defi}

\begin{prop} \label{reflexive-cc}
A $d$-dimensional lattice polytope $P$ containing $0$ in its 
interior is reflexive if and only if 
$F(P) = \{ 0 \}$ and $P$ is canonically closed, i.e., $C(P) = P$. 
\end{prop}

\noindent
{\bf Proof.}  Assume that $P$ is reflexive. Then the interior 
lattice point $0 \in P$ belongs to $F(P)$. On the other hand, 
for every vertex $\nu \in P^*$ one has 
$\ord_P(\nu) = -1$. Therefore the inequalities 
$\langle x, \nu \rangle\geq \ord_P(\nu) +1$ defining $F(P)$ 
are  $\langle x, \nu \rangle\geq 0$  for all vertices $\nu \in P^*$. 
There exists a positive linear combination $\sum_i \lambda_i \nu_i$ of  vertices of $P^*$ 
which is equal to $0$. Therefore, $0 \in M_\R$ is the 
single solution of the inequalities 
\[ \{ x \in M_\R \mid \langle x, \nu \rangle\geq 0 \; \; \forall \nu \in P^* \cap( N \setminus \{0 \}) \}.  \]
It follows that $F(P) = \{ 0 \}$ and $S_F(P)$ contains 
all vertices of $P^*$. By \ref{can-closed-crit}, we obtain 
$C(P) = P$.

Assume now that $F(P) = \{ 0 \}$ and $C(P) = P$. By \ref{can-closed-crit}, all primitive inward pointing  
normals $\nu_Q$ to facets 
$Q \prec P$ belong 
to $S_F(P)$. So we obtain $\ord_P(\nu_Q) =\ord_{F(P)}(\nu_Q) -1 =-1$. Therefore, $P^*$ 
is the convex hull of $S_F(P)$, i.e., $P^*$ is a lattice 
polytope. 
\hfill $\Box$

\begin{rem}
Using the classification due to Kasprzyk \cite{Kas10}, it was shown in  \cite{BKS19} that  there are up to an unimodular isomorphism $661\,280=665\,599 -4\,319$ three-dimensional lattice polytopes  $P$ with 
$F(P) = \{ 0 \}$ which are 
not canonically closed, i.e., they are not reflexive. 
However, the canonical hulls $C(P)$ of these three-dimensional 
polytopes $P$ are always reflexive. 
\end{rem}

\begin{prop} \label{Q-Fine}
Let $P$ be a full dimensional generalized Delzant polytope
with $F(P) \neq \emptyset$ and let $Q \prec P$ be a facet of $P$. 
Assume that $Q$ has nonempty Fine interior, i.e., 
 $F(Q) \neq \emptyset$. 
Then the primitive  inward pointing    
facet normal $\nu_Q$ belongs to $S_F(P)$. In particular, a $d$-dimensional lattice polytope $P$ is canonically closed if  every facet $Q \prec P$  contains at least one  lattice point in its relative interior. 
\end{prop}

\begin{rem} The converse is not true in general. The lattice 
triangle with the vertices $(1,0), (0,1), (-1,-1)$ is canonically closed, but all its facets (sides) have no interior lattice points. 
\end{rem}

\begin{ex}
There are exactly $9$ examples of $3$-dimensional lattice polytopes $P$ with 
$F(P) \neq \emptyset$ and without interior lattice points \cite{BKS19}. All these 
$9$ polytopes are canonically closed, because their facets
contain at least one lattice point in their relative interiors.   
\end{ex}

\noindent
{\bf Proof of \ref{Q-Fine}.} Take a point $q \in F(Q)$ and consider  an arbitrary 
 lattice vector 
$\nu \in S_F(P)$ such that $\nu \not\in \R \nu_Q$. 
Then $\nu$ defines a nonzero element in the lattice $N/\Z \nu_Q$ and we obviously have 
$\ord_Q(\nu) \geq \ord_P(\nu)$. By $q \in F(Q)$, we obtain 
\[ \langle q, \nu \rangle\geq \ord_Q(\nu) +1 \geq \ord_P(\nu) +1. \]
On the other hand, $q$ can not belong to $F(P)$, because 
$q$ is in the relative interior of the facet $Q \prec P$. Hence  
$q$ can not satisfy the inequality  $\langle q, v_Q \rangle 
\geq  \ord_P(\nu_Q) +1$. Therefore, 
the lattice primitive normal $\nu_Q$ must be essential for $F(P)$, 
i.e., $\nu_Q \in S_F(P)$. 
\hfill $\Box$

\begin{coro}
Let $P$ be a full dimensional lattice polytope. Then for any 
integer $k \geq d$ the lattice polytope $kP$ is canonically closed.
\end{coro}

\noindent
{\bf Proof.} The statement follows from \ref{Q-Fine}, because 
for any $(d-1)$-dimensional lattice 
 face  $Q \prec P$  and for any $k \geq \dim Q +1 =d$ the 
lattice polytope $kQ$ has always lattice points in its relative interior.
\hfill $\Box$

\begin{coro}  \label{can-close-f}
Let $P$ be a full dimensional generalized Delzant polytope. 
Then for sufficiently large positive real number $\lambda$ 
the polytope $\lambda P$ is canonically closed. 
\end{coro}

\noindent
{\bf Proof.} By \ref{Q-Fine}, it is enough to choose 
the real number $\lambda \in \R$ such 
that $\lambda Q$ has nonempty Fine interior 
for any facet $Q \prec P$.  
\hfill  $\Box$ 

\begin{prop} \label{can-closed>1}
Let $P$ be a full dimensional generalized Delzant polytope
with $F(P) \neq \emptyset$. If $P$ is canonically closed, then
$\lambda P$ is canonically closed for all $\lambda \geq 1$. 
\end{prop}

\noindent
{\bf Proof.} By \ref{fp-monoton}, $F(\lambda P) \neq \emptyset$
for any $\lambda \geq 1$. 

Consider $\lambda >1$. We claim
that $F(P) + (\lambda -1)P \subseteq F(\lambda P)$. Indeed, take 
$x \in F(P)$ and $y = (\lambda -1)y'$ with $y' \in P$. Then 
for any $\nu \in N \setminus \{0 \}$ we get 
\[ \langle x, \nu \rangle \geq \ord_P(\nu) +1, \;\;  \langle y', \nu \rangle \geq \ord_P(\nu), \;\;  \langle y, \nu \rangle \geq (\lambda -1) \ord_{P}(\nu). \] 
Thus, $ \langle x, \nu \rangle \geq \lambda \ord_P(\nu) +1 = 
 \ord_{\lambda P}(\nu) +1$, 
 i.e. $x + y \in F(\lambda P)$. 
 
 If $\nu \in S_F(P)$ is an inward-pointing primitive lattice 
 normal to a facet $Q \prec P$ such that $\langle x, \nu \rangle =
 \langle y', \nu \rangle +1$ for some $x \in F(P)$ and $y' \in P$, 
 then $\lambda y' \in \lambda P$,  $x + (\lambda -1)y' \in F(\lambda P)$ and 
 \[  \langle x - y',  \nu \rangle =  \langle 
 (x + (\lambda -1)y') - \lambda y' , \nu \rangle = 1,  \]
 i.e., $\nu \in S_F(\lambda P)$. By \ref{can-closed-crit},  $\lambda P$ is canonically closed.  
\hfill  $\Box$ 

\begin{coro}
Let $P$ be a full-dimensional generalized Delzant polytope. Then 
there exists a positive real number $\lambda_P \geq \lambda_{P}^0$ such that 
$\lambda P$ is canonically closed if and only if $\lambda \geq 
\lambda_P$.
\end{coro}

\noindent
{\bf Proof.} Define 
\[ \lambda_P:= \inf \{ \lambda \in \R_{>0} \mid F(P) \neq \emptyset, \;\;  C(\lambda P) = \lambda P \}. \] 
Then $\lambda_P \geq \lambda_P^0$. 

 By \ref{can-closed>1}, 
 $\lambda P$ is canonically closed 
 for any $\lambda > \lambda_P$. 
We claim 
 that also $\lambda_P P$ is canonically closed. The polytopes $\lambda P$ and  $F(\lambda P)$ continuously depend on $\lambda$.  By \ref{can-closed-crit}, we have 
to show that for any primitive inward pointing 
facet normal $\nu_Q$ of $P$ on has  $\nu_Q \in S_F(\lambda_P P)$. Indeed, by taking limit from above,  
we obtain 
\[1 =  \lim_{\lambda \to \lambda_P + 0} \left(\ord_{F(\lambda P)}(\nu_Q) - 
\ord_{\lambda P}(\nu_Q)\right)  =\ord_{F(\lambda_P P)}(\nu_Q) - 
\ord_{\lambda_P P}(\nu_Q) ,   \]
i.e., $\nu_Q \in S_F$. By definition of $\lambda_P$, 
$\lambda P$ is not canonically closed for $\lambda_P^0 \leq \lambda < \lambda_P$.  
\hfill  $\Box$
 
\begin{defi}
Let $P$ be a  $d$-dimensional lattice polytope $P \subset M_\R$ 
with $F(P) \neq \emptyset$. We call 
 $P$  {\em integrally closed} 
if 
$${\rm Conv}(C(P) \cap M) = P.$$ It 
is clear that every canonically closed lattice polytope $P$ 
is integrally closed. 
\end{defi}

\begin{rem}
By \ref{cp-d=2},  every $2$-dimensional lattice polytope  
is simultaneously integrally closed
and canonically closed. 
\end{rem}

\begin{rem}
The $5$-dimensional 
lattice polytope $P$ in Example \ref{can-rat} is integrally closed, but not canonically closed, because $F(P) = \{0 \}$, but 
$P$ is not reflexive (cf. \ref{reflexive-cc}). 
\end{rem}

\begin{rem}
Integrally closed lattice polytopes $P$ with $F(P) = \{0 \}$ are
called in \cite{Bat17} {\em pseudoreflexive}. They satisfy a combinatorial duality that generalize the polar duality for 
reflexive polytopes.  Mavlyutov suggested to use this generalized duality in the Mirror Symmetry for Calabi-Yau complete intersections in toric varieties \cite{Mav11}. 
\end{rem}

\section{The Fine interior and the canonical refinement}

In this section we consider toric varieties associated with noncompact polyhedra \cite[\S 7.1]{CLS11}.  

Let $\sigma \subset N_\R$ be a $d$-dimensional rational finite 
polyhedral cone with the vertex 
$0 =  \sigma \cap -\sigma$.  Consider the $d$-dimensional 
convex lattice polyhedron  
\[ \Theta_\sigma := {\rm Conv}(\sigma \cap (N \setminus  \{0 \} )) \subset \sigma \]
with the recession cone $\sigma$. 
Then the polar dual $d$-dimensional rational polyhedron  
\[ \Theta^\ast_\sigma:= \{ x \in M_\R \mid \langle x, y \rangle  \geq 1 \;\;
\forall y \in \Theta_\sigma \}   \]
is contained in   the dual cone  ${\check \sigma} := 
\{x \in M_\R \mid  \langle x, y \rangle  \geq 0\, \, \forall y \in \sigma \}$ which is the 
recession cone of $\Theta^\ast_\sigma$.  

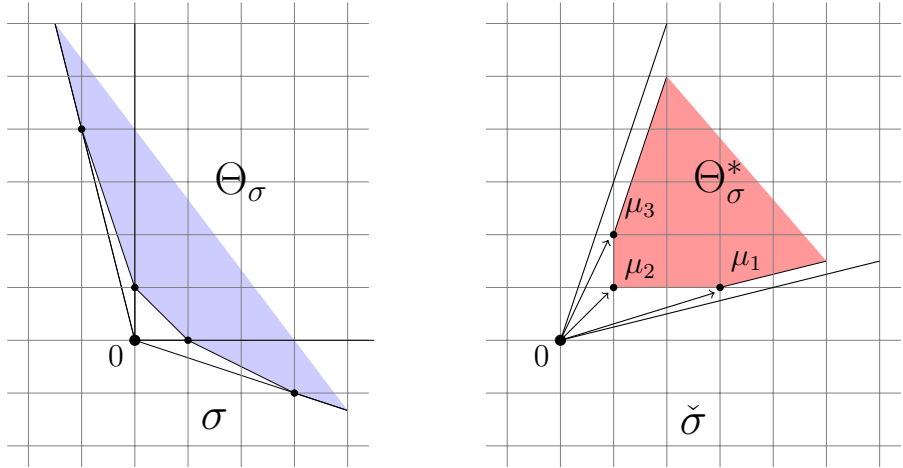
\begin{center}
\begin{tikzpicture}[scale = 0.7]

\begin{scope}[xshift= -2cm, yshift= -2cm] 

\draw[ fill=gray!30!white, very thin] (0,3) -- (-1,0) -- (-1,-1) 
-- ( 1,-1) -- (3,-0.5)   ;

\draw[step=1cm,gray,very thin]
(-3.4,-4.4) grid (4.4,4.4) ;

\fill (-2,-2) circle (3pt);
\fill (-1,-1) circle (2pt);
\fill (-1,0) circle (2pt); 
\fill (1,-1) circle (2pt); 

\draw[ - > ] (-2,-2) -- (-1.1,-1.1) ; 
\draw[ - > ] (-2,-2) -- (0.9,-1.1) ; 
\draw[ - > ] (-2,-2) -- (-1.1,-0.1) ; 

\draw (-2,-2) -- (4, -0.5) ; 

\draw (-2,-2) -- (0, 4);

\node  at (1,1)  {{\Large $\Theta_\sigma^*$}};

\node  at (0.5,-3.5)  {{\Large $\check{\sigma}$}};

\node at (1.5, -0.5) { $\mu_1$ } ; 
\node at (-0.5, -0.66)  { $\mu_2$} ;
\node at (-0.5, 0.5)  { $\mu_3$}; 

\node[left]  at (-2,-2.3) {$0$};

\end{scope} 
\begin{scope}[xshift= -12cm ] 

\draw[ fill=gray!20!white, very thin] (-1.5,2) -- (-1,0) -- (0,-3) 
-- ( 1,-4) -- (3,-5) -- (4,-5.33)  ;

\draw[very thin] (0,-4) -- (4,-5.33); 
\draw[very thin] (0,-4) -- (1,-4);  
\draw[very thin] (0,-4) -- (0,-3);
\draw[very thin] (0,-4) -- (-1,0);
\draw[very thin] (0,-4) -- (-1.5, 2);

\draw[step=1cm, gray, very thin]
(-2.4,-6.4) grid (4.4,2.4) ;

\draw[very thin] (0,-4) -- (4.5,-4); 
\draw[very thin] (0,-4) -- (0,2); 

\fill (0,-4) circle (3pt);
\fill (1,-4) circle (2pt);
\fill (0,-3) circle (2pt);
\fill (3,-5) circle (2pt);
\fill (-1,0) circle (2pt); 

\node [left] at (0,-4.3) {{0}}; 

\node  at (2,-1)  {{\Large $\Theta_\sigma$}};

\node  at (1.5,-5.5)  {{\Large ${{\sigma}}$}};

\end{scope}
\end{tikzpicture}
\end{center}

\noindent
One has a natural bijection between 
compact facets $\theta_i \prec \Theta_\sigma$ 
and rational 
vertices  $\mu_i \in \Theta^\ast_\sigma$ such that 
\[ \theta_i = \{ y \in \Theta_\sigma \, |\, 
\langle \mu_i, y \rangle = 1 \}. \]
Denote by  ${\mathcal M}({\sigma}) := \{ \mu_1, \ldots, \mu_s \}$  the set of all (rational) vertices of the polyhedron $\Theta^\ast_\sigma$. 
Then ${\mathcal M}({\sigma})$ belongs to the interior of the 
recession cone $\check{\sigma}$  and we have 
\[ \Theta_\sigma = \{ y \in \sigma  \,  | \, 
\langle \mu , y 
\rangle \geq 1, \;\; \forall \mu_i \in  {\mathcal M}({\sigma})\}, \;   \;\; \Theta^\ast_\sigma   ={\rm Conv}\left({\mathcal M}(\sigma) \right) + \check{\sigma} . 
\] 

\begin{defi} 
Consider  the upper-convex piecewise linear 
function on $\sigma$:
\[ \gamma_\sigma(y) := \ord_{\Theta^\ast_\sigma}(y) = \min_{x \in \Theta^\ast_\sigma}  \langle x, y 
\rangle = 
\min_{\mu_i \in {\mathcal M}({\sigma}) }  \langle \mu_i , y 
\rangle, \;\; \forall y \in \sigma. \] 
The domains of linearity of the function 
$\gamma_\sigma$ 1-to-1 correspond to 
rational points $\mu_i \in {\mathcal M}({\sigma})$ and they 
define a fan ${\sigma}_{\rm can}$ 
which is called the 
{\em   canonical refinement} of the full dimensional cone 
$\sigma$. The 
canonical refinement 
 ${\sigma}_{\rm can}$ is the normal fan of the convex 
 polyhedron ${\Theta^\ast_\sigma}$ \cite[Def. 7.1.3]{CLS11}.
\end{defi}
 
\begin{prop} \cite[Prop. 11.4.15]{CLS11}
Let $X_\sigma$ be the affine toric
variety of a $d$-dimensional   
cone $\sigma \subset N_\R$. Then  the canonical refinement 
${\sigma}_{\rm can}$ of $\sigma$ defines a quasi-projective 
toric variety $X_{{\sigma}_{\rm can}}$ 
together with a proper toric morphism 
\[ \phi_\sigma\, :\, 
X_{{\sigma}_{\rm can}} \to X_\sigma \]
such that the toric variety $X_{{\sigma}_{\rm can}}$ has 
$\Q$-Gorenstein canonical singularities and the canonical 
class $K_{X_{{\sigma}_{\rm can}}}$ is an ample $\Q$-Cartier 
divisor defined by the upper-convex ${\sigma_{\rm can}}$-piecewise linear function 
 $\gamma_\sigma = \ord_{\Theta^\ast_\sigma}$. In particular, 
 $\Theta^\ast_\sigma$ is the supporting 
 polyhedron of the canonical class on $X_{{\sigma_{\rm can}}}$. 
\end{prop}

Our interest to the canonical refinements is motivated by 
the following statement: 

\begin{thm} \label{FP-refinement}
Let $\Sigma_P$ be the normal fan of  a d-dimensionsional 
generalized Delzant polytope $P \subset M_\R$, i.e., 
\[ P = \bigcap_{\sigma \in \Sigma_P(d)} 
(p_\sigma + \check{\sigma} ), \]
where $p_\sigma \in P$ denotes the vertex of $P$ corresponding 
to a $d$-dimensional cone $\sigma \in \Sigma_P(d)$.   
Then 
\[ F(P) = \bigcap_{\sigma \in \Sigma_P(d)} 
(p_\sigma + \Theta^\ast_{\sigma} ). \]   
\end{thm}
 
\noindent
{\bf Proof.} Take a $d$-dimensional cone 
$\sigma \in \Sigma_P(d)$. Then for any 
nonzero lattice vector $\nu \in \sigma\cap N$ we have 
$\ord_P(\nu) = 
\langle p_\sigma, \nu \rangle$, where $p_\sigma$ is the 
vertex of $P$ corresponding to $\sigma$. Hence,  
we obtain  
\[ p_\sigma + \Theta^\ast_\sigma = 
\{ x \in M_\R \mid  
\langle x, \nu \rangle \geq \ord_P(\nu) +1 \;\; \forall 
\nu \in \sigma \cap N \setminus \{0 \} \},  \]
because $\langle x, \nu \rangle \geq \ord_P(\nu) +1$
is equivalent to $\langle x - p_\sigma, \nu \rangle \geq 1$. Since $\Sigma_P$ is a complete fan, every 
lattice vector $\nu \in N$ is contained in some 
$d$-dimensional cone $\sigma$. This implies the 
statement. 
\hfill $\Box$

\begin{coro} \label{gen-fp}
Let $\sigma \in \Sigma_P[d]$ be a full dimensional 
cone of the normal fan $\Sigma_P$. We denote by 
${\sigma}_{\rm can}[1]$ be the set of all 
primitive lattice generators of $1$-dimensional cones 
in the fan ${\sigma}_{\rm can}$. Then 
\[ F(P) = \{ x \in M_\R \mid 
\langle x, \nu \rangle \geq \ord_P(\nu) +1 \;\; \forall 
\nu \in  \bigcup_{\sigma \in \Sigma_P(d)} {\sigma}_{\rm can}[1] \}.  \]
\end{coro}

One can extend the notion of the canonical refinement 
from "local" to "global" by gluing  
the canonical refinements $\sigma_{\rm can}$ of $d$-dimensional cones $\sigma 
\in \Sigma(d)$ of  
a complete fan $\Sigma$. Using the equalities  
\[ \Theta_{\sigma} \cap \Theta_{\sigma'} = 
\Theta_{\sigma \cap \sigma'}, \;\; \forall \sigma, \sigma' \in \Sigma(d),  \] 
we can glue local canonical refinements $\sigma_{\rm can}$ 
 of cones $\sigma \in \Sigma(d)$  
and obtain a global 
canonical refinement ${\Sigma}^{\rm can}$ 
of the complete fan $\Sigma$. 
\medskip

In the next section we will prove the following  
general statement 
about the canonical refinement $\Sigma_P^{\rm can}$ 
of the normal fan $\Sigma_P$ of a canonically closed full dimensional generalized 
Delzant polytopes $P$:

\begin{thm} \label{Mink-sum} 
Let 
$P$ be a canonically closed 
full dimensional generalized Delzant polytope with 
$F(P) \neq \emptyset$. Then 
for any $\lambda > 1$ the normal fan $\Sigma_{F(\lambda P)}$ 
of the 
full dimensional polytope $F(\lambda P)$ is the canonical 
refinement $\Sigma^{\rm can}_P$ 
of the normal fan $\Sigma_P$ of 
$P$ and 
one has the Minkowski sum decomposition
\[ F(\lambda P) = F(P) + (\lambda -1) P. \]
\end{thm}

\begin{coro} \label{Mink-sum2}
Let 
$P$ be a canonically closed 
full dimensional generalized Delzant polytope with 
$F(P) \neq \emptyset$. Then 
one has the Minkowski sum decomposition
\[ F(2 P) = F(P) + P. \]
\end{coro}

\noindent
{\bf Proof.} The statement follows from  Theorem 
\ref{Mink-sum} if one takes  $\lambda = 2$. \hfill $\Box$

\begin{coro}
Let $P$ be an arbitrary $d$-dimensional generalized 
Delzant polytope and let $\Sigma^{\rm can}_P$ be the canonical refinement  of its normal fan $\Sigma_P$. Then 
$\Sigma^{\rm can}_P$ is the normal fan  $\Sigma_{F(\lambda P)}$ 
of the Fine interior $F(\lambda P)$  $\forall \lambda \gg  1$.
\end{coro}

\noindent
{\bf Proof.} For  $\lambda \gg 1$ the polytope  $F(\lambda P)$ is full dimensional and  the polytope $\lambda P$ is 
canonically closed (see \ref{can-close-f} and \ref{can-closed>1}). \hfill $\Box$

\begin{prop} \label{gg1}
Let $P$ be a canonically closed $d$-dimensional generalized 
Delzant polytope with $F(P) \neq \emptyset$. Then for any 
$\lambda > 1$  the combinatorial type of the polytope $F(\lambda P)$ 
is independent on 
$\lambda$ and the set of vertices of $F(\lambda P)$ equals 
\[ \bigsqcup_{\sigma \in \Sigma_P(d)} 
\{ \lambda p_\sigma + \mu_i \; \mid \; \mu_i \in \mathcal M(\sigma)\}. \]
\end{prop}

\noindent
{\bf Proof.} By Theorem \ref{Mink-sum}, for  $\lambda > 1$  the normal fan
of $F(\lambda P)$ is just the canonical refinement of $\Sigma_{P}$ which is independent on  $\lambda > 1$. The combinatorial structure  of the  polytope $F(\lambda)$ is dual to the one of its normal fan. So it is also independent on $\lambda >1$. Using the Minkowski sum decomposition 
\[ F(\lambda P) = (\lambda -1) P + F(P), \;\; \lambda >1, \]
we see that every vertex $q \in F(\lambda P)$ has a unique 
representation as sum 
\[ q = (\lambda -1)p_\sigma + r = \lambda p_\sigma + (r - p_{\sigma}), \]
where  
$p_\sigma$ is a vertex of $P$ corresponding to some 
$d$-dimensional cone $\sigma \in \Sigma_P(d)$ and $r$ is 
a vertex of $F(P)$ such that $r- p_{\sigma} \in \check{\sigma}$. Note that  the vertex $q$ of $F(\lambda P)$ corresponds to a full dimensional cone of 
the canonical refinement $\Sigma_P^{\rm can}$ which 
is contained in $\sigma$ and spanned by some compact facet $\theta_i \prec \Theta_\sigma$. 
Therefore the vertex $r \in F(P) \subset P$ determines the facet $\theta_i \prec \Theta$ by the condition 
\[ \theta_i = \{ y \in \Theta_\sigma  
\; \mid \; \langle r- p_{\sigma}, y \rangle =1 \},    \]  
i.e., $r - p_\sigma = \mu_i \in {\mathcal M}(\sigma)$ and $q = \lambda p_\sigma + \mu_i$.    
\hfill $\Box$ 

\begin{coro}
Let $P$ be any $2$-dimensional lattice polytope with 
$F(P) \neq \emptyset$. Then 
$F(P)$ is the convex hull of the lattice points  
$$\bigcup_{\sigma \in \Sigma_P(2)}  \{ p_\sigma + \mu_i\; |  
\; \mu_i \in {\mathcal M}({\sigma})\} .$$
\end{coro}

\noindent
{\bf Proof.} We know that if $P$ is a $2$-dimensional lattice polytope with $F(P) \neq \emptyset$ then $P$ is canonically closed and  
$F(P) = {\rm Conv}(Int(P) \cap M)$. For any $2$-dimensional cone $\sigma \in \Sigma_P(2)$ the 
set ${\mathcal M}(\sigma)$ consists of lattice vectors, because the canonical refinement $\Sigma_P^{\rm can}$ of $\Sigma_P$ defines 
a $2$-dimensional Gorenstein toric variety. By \ref{gg1}, for all  $\lambda >1$ the set of vertices of the Fine interior $F(\lambda P)$ equals 
\[ \bigsqcup_{\sigma \in \Sigma_P(2)} 
\{ \lambda p_\sigma + \mu_i \; \mid \; \mu_i \in \mathcal M(\sigma)\}. \]
Taking the limit $\lambda \to 1$, we obtain that all vertices of 
$F(P)$ are contained in the finite set of the limiting  lattice points 
\[ \bigcup_{\sigma \in \Sigma_P(2)} 
\{ p_\sigma + \mu_i \; \mid \; \mu_i \in \mathcal M(\sigma) \}  \subset Int(P) \cap M. \]
Note that after taking the limit $t \to 1$ two disjoint finite subsets 
\[ \{ p_\sigma + \mu_i \; \mid \; \mu \in \mathcal M(\sigma)\}, \;\;  \{ p_{\sigma'} + \mu_j \; \mid \; \mu_j \in \mathcal M(\sigma')\} \]
corresponding to different $2$-dimensional cones 
$\sigma, \sigma' \in \Sigma_P(2)$ may get common 
elements. The latter happens, i.e.,  for all $2$-dimensional reflexive polygons.  
\hfill $\Box$ 

\section{The canonical toric variety $\widetilde{V}$}

\begin{defi}
Let $P$ be a full dimensional generalized Delzant 
polytope with $F(P) \neq \emptyset$.  Consider the continuous 
real  function 
$\delta_P\, :\, N_\R \to \R$,  
$$\delta_P(y):= 
\ord_{F(P)}(y)  - \ord_P(y), \;\;  y \in N_\R.$$  
\end{defi}

\begin{prop} \label{alpha-a-e}
The function $\delta_P$ is linear on cones of the normal fan 
of the Minkowski sum $P + F(P)$ and 
it has the following properties: 

{\rm (a)}  $\delta_P(y) \geq 0$  for all $y \in N_\R;$

{\rm (b)}  $\delta_P(y) = 0$  if and only if  $ y =0;$

{\rm (c)} $\delta_P(\lambda y) = \lambda \delta_P(y)$ for 
all $\lambda  \in \R_{\geq 0};$ 

{\rm (d)}  $ \{ \nu \in N \;  | \; \delta_P(\nu) = 1 \} = S_F(P);$

{\rm (e)} $\{ \nu \in N \; |  \; 0 < \delta_P(\nu) < 1 \} = \emptyset$.

In particular,  
\[ \Delta_P:= \{ y \in N_\R \; | \; \delta_P(y) \leq 1 \} \]
is a full dimensional compact subset in $N_\R$ containing  $S_F(P)$ on its boundary.
\end{prop}

\noindent
{\bf Proof.} The normal fan $\Sigma'$ 
of the Minkowski sum $P + F(P)$ 
is the coarsest common refinement of the normal fans 
$\Sigma_P$ and 
$\Sigma_{F(P)}$ \cite[Prop. 6.2.13, \S 6]{CLS11}. Therefore, 
the upper convex functions $\ord_P$ and $\ord_{F(P)}$ are linear on each cone of the normal fan 
$\Sigma'$.  This implies (c). Since the polytope $F(P)$ is strictly contained in the interior 
of $P$, we have $\ord_{F(P)}(y) > \ord_P(y)$ for all $y \in 
N_\R \setminus \{0 \}$. This implies (a) and (b). The statement 
(d) follows from the definition of $S_F(P)$, and  (e) follows from 
the inequality $\ord_{F(P)}(\nu) - \ord_P(\nu) \geq 1$ for all nonzero 
$\nu \in N$. 
\hfill $\Box$

We note that the compact set  $\Delta_P$ is usually not convex (see pictures in  Corollary \ref{2-dim-can} and Example \ref{enri}).

\begin{thm} \label{can-sing1}
Let $P$ be a full dimensional generalized Delzant polytope $(F(P) \neq \emptyset)$. Denote by $\widetilde{\Sigma}$ 
the normal fan of the Minkowski 
sum $\widetilde{P} := F(P) + C(P)$. 
Then any primitive inward pointing facet normal $\nu \in  \widetilde{\Sigma}[1]$ is contained in the support 
of the Fine interior $S_F(P)$. 
\end{thm}

\noindent
{\bf Proof.} Since all generalized Delzant  polytopes $\lambda P$
($\lambda \in \R_{>0})$ have the same normal fan, by 
\ref{Fp-computer}, there exists 
a finite subset $S \subset N$ such that  $S_F(\lambda P) \subseteq S$ for all $ \lambda \geq 1$. Recall that 
for any $\lambda >1$ the Fine interior $F(\lambda P)$ is full dimensional, because it contains the full dimensional 
Minkowski sum 
$F(P) + (\lambda -1)P$ (see the proof of \ref{can-closed>1}). 
On the other hand,  any primitive inward pointing facet 
normal $\nu_Q$ of the polytope $F(\lambda P)$
belongs to $S_F(\lambda P)$,  i.e.,  $\nu_Q \in S$. 
For sufficiently small $\varepsilon >0$ the set 
$\{ \nu_1, \ldots, \nu_k \} = \Sigma[1]:= \Sigma_{F((1 + \varepsilon) P)}[1] \subset N$ of all  primitive
inward pointing facet normals to 
 facets $Q_1(\varepsilon), \ldots, Q_k(\varepsilon)$ 
of the full dimensional polytope 
$F((1 + \varepsilon)P)$ does not depend on $\varepsilon$ and 
the corresponding supporting affine hyperplanes 
$H_i(\varepsilon)$ spanned by facets $Q_i(\varepsilon)$ 
contribute to the Fine interior $F((1 + \varepsilon)P)$ 
of the polytope
 $(1+ \varepsilon)P$, i.e., $\Sigma[1]= \{ \nu_1, \ldots, \nu_k \} 
\subseteq  
S_F((1+ \varepsilon)P)$ and 
\[ \ord_{F((1+ \varepsilon)P)}(\nu_i) -  
\ord_{(1+ \varepsilon)P}(\nu_i) = 1  \; \; \forall i \in \{ 1,\ldots, k \}. \]
The polytopes $(1+ \varepsilon)P$ and $F((1 + \varepsilon)P)$ depend continuously on $\varepsilon$. Hence,  by taking the limit $\varepsilon \to 0$,
we obtain that 
\[ \ord_{F(P)}(\nu_i) -  
\ord_{P}(\nu_i) = 1  \; \; \forall  i \in \{ 1,\ldots, k \}, \]
i.e., $\{ \nu_1, \ldots, \nu_k \} =  
\Sigma[1]\subseteq S_F(P)$.  

Next we want to show that for sufficiently 
small $\varepsilon >0$ one has the following Minkowski sum 
decomposition 
\[  F((1 + \varepsilon)P) = F(P) + \varepsilon C(P). \]
Note that the full dimensional polytope   $F((1 + \varepsilon)P)$ 
is the intersection of the $k$ half-spaces defined by primitive 
lattice facet normals $\nu_1, \ldots, \nu_k$. So we have 
\begin{align*}
 F((1 + \varepsilon)P) = & \{  x \in M_\R \; | \; \langle x, \nu \rangle 
\geq \ord_{(1+ \varepsilon)P}(\nu) + 1, \;\; 
\forall \nu \in \Sigma[1] \} = \\  
 = & \{  x \in M_\R \; | \; \langle x, \nu \rangle 
\geq (1+ \varepsilon) \ord_{P}(\nu) + 1, \;\; 
\forall \nu \in \Sigma[1] \} =  \\
= & \{  x \in M_\R \; | \; \langle x, \nu \rangle 
\geq (\ord_{P}(\nu) + 1) + \varepsilon \ord_{P}(\nu), \;\; 
\forall \nu \in \Sigma[1] \}. 
\end{align*}
This implies that 
\[ F(P) +\varepsilon C(P) \subseteq F((1 +\varepsilon)P \]
and  that the full dimensional polytope $F((1 + \varepsilon)P)$  
defines an ample $\R$-Cartier divisor 
\[ D(\varepsilon):= \sum_{i=1}^k - (1 + \varepsilon) \ord_P(\nu_i) D_i - \sum_{i=1}^k D_i   \]
 on the projective toric 
variety $V:= V(\Sigma)$ of the  fan 
$\Sigma$ which is independent 
of the small $\varepsilon >0$. In particular, 
$D(\varepsilon) - D(\varepsilon') = 
(\varepsilon'- \varepsilon)\sum_{i=1}^k  \ord_P(\nu_i) D_i$ 
is an $\R$-Cartier divisor on $V$. Therefore, $$D:= -\sum_{i=1}^k  \ord_P(\nu_i) D_i$$ is an $\R$-Cartier divisor on $V$. 

We claim that any primitive inward pointing facet normal $\nu 
\in S_F(P)$ of $C(P)$ 
is contained in $\Sigma[1]$. Indeed, let $\Gamma_\nu \subset C(P)$
be a facet of $C(P)$ defined by the equation 
$\langle x, \nu \rangle = \ord_P(\nu)$.
Then there exist a point $x' \in F(P)$ such that 
$\langle x', \nu \rangle = \ord_P(\nu) +1$. 
The  $(d-1)$-dimensional 
polytope $x' + \varepsilon \Gamma_\nu$ is contained 
in $F(P) + \varepsilon C(P) \subseteq F((1+\varepsilon)P)$ 
and $$\langle x'' , \nu \rangle  = (1 + \varepsilon )\ord_P +1  = 
\ord_{(1 + \varepsilon)P}(\nu) + 1 \;\; \forall 
x'' \in x' + \varepsilon \Gamma_\nu. $$
Hence, $ x' + \varepsilon \Gamma$ is a facet of $ F((1+\varepsilon)P)$ and $\nu \in \Sigma[1]$. 

Since the full dimensional polytope $C(P)$ equals the intersection 
of all half-spaces $\langle *, \nu \rangle \geq \ord_P(\nu)$ defined 
by 
primitive inward pointing facet normals  $\nu 
\in S_F(P)$,  and $\ord_{C(P)}(\nu) = \ord_P(\nu)$ $\forall \nu \in 
S_F(P)$ (see \ref{can-ord} (a)), we obtain
\[ C(P) =  \{ x \in M_\R \mid \langle x, \nu \rangle \geq 
\ord_P(\nu), \; \; 
\forall \nu \in \Sigma[1] \}. \]
i.e.,  that the $\Sigma$-piecewise 
linear function corresponding to the $\R$-Cartier 
divisor $D$  equals the  upper-convex $\Sigma$-piecewise 
linear function function 
$\ord_{C(P)}$. 

On the
other hand, $\ord_{F(P)}$ is the limit $\varepsilon \to 0$ of  convex
$\Sigma$-piecewise linear functions 
$\ord_{F((1+ \varepsilon)P)}$. Therefore, $\ord_{F(P)}$
is also an  upper-convex $\Sigma$-piecewise linear function. 
Now we obtain that for sufficiently small  
$\varepsilon >0$ the strictly upper-convex $\Sigma$-piecewise 
linear function $\ord_{F((1+ \varepsilon)P)}$ equals the sum of 
two upper-convex  $\Sigma$-piecewise 
linear functions:  
\[ \ord_{F((1+ \varepsilon)P)} = 
\ord_{F(P)} +  \varepsilon \ord_{C(P)}. \]
Therefore,  $F((1 + \varepsilon)P) = F(P) + \varepsilon C(P)$ 
and $\Sigma$ is the coarsest common refinement of the normal
fans $\Sigma_{F(P)}$ and $\Sigma_{C(P)}$
(see \cite[Prop. 6.2.13]{CLS11}). In particular, 
the projective fan $\Sigma$ with $\Sigma[1] \subset S_F(P)$ equals the normal fan $\widetilde{\Sigma}$ of the 
Minkowski sum $\widetilde{P} := C(P) + F(P)$, i.e., $\widetilde{\Sigma}[1] \subset S_F(P)$. \hfill  $\Box$

\begin{coro} \label{normal-fan_tilde}
Let $P \subset M_\R$ be a full dimensional lattice polytope
with $F(P) \neq \emptyset$ and let $\widetilde{P} := F(P) + C(P)$. 
Then the normal fan of $\widetilde{P}$ equals the normal 
fan of $F((1 + \varepsilon)P)$ for sufficiently small positive 
$\varepsilon$. 
\end{coro}

\noindent
{\bf Proof.} The statement immediately follows from the proof of Theorem \ref{can-sing1}. \hfill $\Box$

\begin{coro} \label{can2}
Let $P \subset M_\R$ be a full dimensional lattice polytope
with $F(P) \neq \emptyset$ and let $\widetilde{P} := F(P) + C(P)$.  Denote by $\widetilde{V}$ the projective 
toric variety of 
the normal fan $\widetilde{\Sigma}$ of $\widetilde{P}$. Then 
$\widetilde{V}$ is a $\Q$-Gorenstein 
toric variety with at worst canonical singularities. 
\end{coro}

\noindent
{\bf Proof.}
Since $\widetilde{\Sigma}$ is a common refinement 
of the normal fans $\Sigma_{C(P)}$ and $\Sigma_{F(P)}$,  
both functions 
$\ord_{C(P)}$ and $\ord_{F(P)}$ are $\widetilde{\Sigma}$-piecewise linear. It follows that  
$$\delta_{C(P)}:= 
\ord_{F(P)} - \ord_{C(P)}$$ is also $\widetilde{\Sigma}$-piecewise
linear. By \ref{can-ord}(a), $\widetilde{\Sigma}$-piecewise
linear function 
$\delta_{C(P)}$ has value 
$1$ on every primitive lattice vector 
$\nu \in \widetilde{\Sigma}[1] \subseteq S_F(P)$. Therefore, 
the canonical class of  $\widetilde{V}$ is a 
$\Q$-Cartier divisor defined by 
the $\widetilde{\Sigma}$-piecewise linear 
function $\delta_{C(P)}$. By \ref{can-ord} and \ref{can-suppF}, 
we can apply \ref{alpha-a-e} to the rational 
polytope $C(P)$  and 
obtain 
\[ \delta_{C(P)}(\nu) = 
\ord_{F(P)}(\nu) - \ord_{C(P)}(\nu) \geq 1, \;\; \forall \nu \in N \setminus \{0\}. \] 
Therefore, all singularities of the toric variety $\widetilde{V}$ are 
at worst canonical. 
\hfill $\Box$

\begin{coro} \label{can3-crep}
Let $\widehat{\Sigma}$ be a projective refinement 
of the fan $\widetilde{\Sigma}$ such that 
$\widehat{\Sigma}[1] = S_F(P)$. Then the corresponding 
projective 
toric morphism $\phi\, :\, \widehat{V} \to \widetilde{V}$ 
is crepant. Furthermore, if $\psi \, :\, V \to \widetilde{V}$ is  a toric desingularization 
of $\widetilde{V}$ defined by a regular refinement 
$\Sigma$ of $\widetilde{\Sigma}$,  then 
\[ K_V = \psi^\ast K_{\widetilde{V}} + \sum_{\nu \in \Sigma[1] \setminus \widetilde{\Sigma}[1]} \widetilde{a}(\nu) D_\nu 
,  \]
where $ \widetilde{a}(\nu) := \delta_{C(P)}(\nu) -1 = \ord_{F(P)}(\nu)  - \ord_{C(P)}(\nu) -1$.  In particular, 
 $\phi\, :\, \widehat{V} \to \widetilde{V}$ is
a maximal projective partial crepant desingularization of $\widetilde{V}$. 
\end{coro}

\noindent
{\bf Proof.}
By \ref{can-ord}(a), the $\widetilde{\Sigma}$-piecewise linear 
function $\delta_{C(P)} := \ord_{F(P)} - \ord_{C(P)}$ representing 
the canonical class $K_{\widetilde{V}}$ takes value 
$1$ on every lattice vector $\nu \in S_F(P) = 
\widehat{\Sigma}[1]$. Therefore, $\phi^*K_{\widetilde{V}} = 
K_{\widehat{V}}$, i.e., $\varphi$ is crepant. 

If $\psi \, :\, V \to \widetilde{V}$ is  an arbitrary 
toric desingularization 
of $\widetilde{V}$ defined by a regular refinement 
$\Sigma$ of $\widetilde{\Sigma}$, then the discrepancy  $\widetilde{a}(\nu)$ for 
a divisorial valuation $\nu \in N$ is the difference between 
the values of the $\Sigma$-piecewise linear functions 
corresponding to $K_V$ and $K_{\widetilde{V}}$, i.e., 
\[ \widetilde{a}(\nu)  = \delta_{C(P)}(\nu) -1  = (\ord_{F(P)}(\nu)  - \ord_{C(P)}(\nu)) -1. \]
In particular,  $\widetilde{a}(\nu) = 0$ if and only if 
$\nu \in S_F(P)$, i.e., $\phi$ is
a maximal projective partial crepant desingularization of $\widetilde{V}$. 
\hfill $\Box$ 
\medskip

\begin{thm} \label{codim2faces} Let $P \subset M_\R$ be a full dimensional lattice polytope
with $F(P) \neq \emptyset$, $\widetilde{P} = C(P) + F(P)$, 
and $\widetilde{\Sigma}$ the normal fan of $\widetilde{P}$.
Consider an arbitrary $2$-dimensional cone $\sigma \in \widetilde{\Sigma}$ and denote by  $\nu_i, \nu_j \in \widetilde{\Sigma}[1] \subseteq S_F(P)$ its spanning primitive lattice vectors. Then \\

{\rm (a)} the $(d-2)$-dimensional affine subspace $L_\sigma^1 \subset M_\R$ defined by the equations 
\[ \langle x, \nu_i \rangle = \ord_P(\nu_i) +1, \;\;  
\langle x, \nu_j \rangle = \ord_P(\nu_j) +1 \]
has nonempty intersection with $F(P)$,

{\rm (b)} 
 the $(d-2)$-dimensional affine 
subspace $L_\sigma^0 \subset M_\R$ defined by the equations 
\[ \langle x, \nu_i \rangle = \ord_P(\nu_i), \;\;  
\langle x, \nu_j \rangle = \ord_P(\nu_j) \]
contains at least one vertex of $P$.
\end{thm}

\noindent{\bf Proof.}  (a) By \ref{normal-fan_tilde}, the normal fan of $\widetilde{P}$ equals 
the normal fan of the Fine interior  of $(1 + \varepsilon)P$ 
for small $\varepsilon >0$. Therefore, for these $\varepsilon$, the $d$-dimensional 
polytope $F((1 + \varepsilon)P)$
has a $(d-2)$-dimensional face $\Gamma_{\sigma, \varepsilon}$ whose 
affine span is the intersection $L_{\sigma, \varepsilon}^1$ of two affine hyperplanes defined by the equations
\[  \langle x, \nu_i \rangle = \ord_{(1+\varepsilon)P}(\nu_i) + 1, \;\;  \langle x, \nu_j \rangle = \ord_{(1+\varepsilon)P}(\nu_j) + 1.\]
Assume that the intersection $L^1_\sigma \cap F(P)$ is empty. 
This means that the distance between the compact convex 
set $F(P)$ and the affine linear subspace $L^1_\sigma$ is positive.
On the other hand, we have 
\[ L_{\sigma}^1 = \lim_{\varepsilon \to 0} L_{\sigma, \varepsilon}^1, \;\;
F(P) =\lim_{\varepsilon \to 0} F((1 + \varepsilon)P). \]
This imply that for sufficiently small $\varepsilon >0$ the 
intersection $  L_{\sigma, \varepsilon} \cap F((1 + \varepsilon)P)$ must be also empty. The latter  contradicts the fact that 
 $(d-2)$-dimensional face $\Gamma_{\sigma, \varepsilon}$ of $F((1 + \varepsilon)P)$ 
  is contained in $L_{\sigma, \varepsilon}^1$ for all sufficiently small 
$\varepsilon > 0$.  

(b) Consider 
the $2$-dimensional  sublattice $N_\sigma \subset N$ spanned 
by $N \cap \sigma$. Then $N_\sigma \cong \Z^2$ is a direct summand of 
$N$. Let  $e_1, e_2$ be a $\Z$-basis of $N_\sigma$.  Consider  the lattice projection 
$$\pi_\sigma \, :\, M \to \Z^2, \; \; 
m \mapsto (\langle m, e_1 \rangle,  \langle m, e_2 \rangle).$$ 
Since $\widetilde{V}$ has at worst canonical singularities, the lattice triangle with vertices $0, \nu_i, \nu_j$ has no interior 
lattice points,  and one can choose a $\Z$-basis $e_1, e_2 \in N_\sigma$ in  
such a way that $\nu_i =e_1$ and $\nu_j = e_1 + k e_2$ for some 
positive integer $k$. We claim that all $k+1$ lattice vectors
\[ e_1,  e_1 + e_2, e_1 + 2e_2, \ldots, e_1 + (k-1) e_2, e_1 + ke_2  \]
belong to $S_F(P)$. For the lattice vectors $e_1 = \nu_i$ and $e_1 + ke_2 = \nu_j$ this follows from \ref{can-sing1}. By \ref{can-ord}(a), we have
\[   {\ord}_{C(P)} (\nu_i) +1 = {\ord}_P (\nu_i) +1 = {\ord}_{F(P)}(\nu_i),\]  
\[   {\ord}_{C(P)} (\nu_j) +1 =
{\ord}_P (\nu_j) +1 = {\ord}_{F(P)}(\nu_j).  \]
The functions ${\ord}_{F(P)}$ and ${\ord}_{C(P)}$ are linear on the $2$-dimensional cone $\sigma$. This implies 
\[  {\ord}_{C(P)} (\nu) +1 = {\ord}_{F(P)}(\nu) \] 
for any lattice point $\nu$ of the segment $[\nu_i, \nu_j]$. 
By \ref{can-suppF}, all these points belong to $S_F(P) = S_F(C(P))$. In particular, 
we have 
\[  {\ord}_{P} (\nu) +1 = {\ord}_{F(P)}(\nu), \;\; \forall \nu \in 
N_\sigma \cap [\nu_i, \nu_j], \]
i.e., $e_1 + re_2 \in S_F(P)$ for all $r \in \{ 0, 1, \ldots, k \}$. 

Now apply (a) and take a point $p_0 \in F(P) \cap L_\sigma^1$.  Then 
\[ \langle p_0, \nu_i \rangle  = {\ord}_P(\nu_i) + 1 ={\ord}_{F(P)}(\nu_i) ,  \;\; 
\langle p_0, \nu_j \rangle  = {\ord}_P(\nu_j) + 1 ={\ord}_{F(P)}(\nu_j).\]

By the linearity of ${\ord}_{F(P)}$ on $\sigma$, we obtain 
\[ \langle p_0, \nu \rangle  = {\ord}_P(\nu) + 1 \;\; \forall \nu \in 
N_\sigma \cap [\nu_i, \nu_j].  \] 
In particular, this implies  
\[  \langle p_0, \nu \rangle \in \Z \;\; \forall \nu \in 
N_\sigma \cap [\nu_i, \nu_j]. \]

Since the set $N_\sigma \cap [\nu_i, \nu_j]$ generate the sublattice $N_\sigma \subseteq N$, we obtain that the $\pi_\sigma$-projection of $p_0 \in M_\Q$ to $\Q^2$ is a lattice point $\overline{p_0} \in \Z^2$. 

Obviously, the projection $P_\sigma:= \pi_\sigma (P) \subset \R^2$
is a lattice polygon contained in the angle $A_\sigma$ defined by 
two inequalities $\langle x, \nu_i \rangle  \geq \ord_P(\nu_i)$ and $\langle x, \nu_j \rangle  \geq \ord_P(\nu_j)$:
\[ P_\sigma \subset A_\sigma =\{ (x_1, x_2) \in \R^2 \, \mid \, x_1 \geq \ord_P(\nu_i),\; \; x_1 + k x_2 \geq \ord_P(\nu_j) \}, \]
and both sides of $A_\sigma$ contain lattice vertices of $P_\sigma$. 
Moreover, the lattice point $\overline{p_0} \in \Z^2 \subset \R^2$ is the intersection of two lines with the equations
\[ x_1 = \ord_P(\nu_i) +1, \;\; x_1 + kx_2 = \ord_P(\nu_j) +1. \]
Therefore, the intersection of two another lines 
with the equations 
\[ x_1 = \ord_P(\nu_i), \;\; x_1 + kx_2 = \ord_P(\nu_j)\]
is also a lattice point $q \in \Z^2$ which is the vertex 
of the angle $A_\sigma \subset \R^2$. By shifting 
the lattice polytope $P$, we can now assume without loss 
of generality that $\ord_P(\nu_i) = \ord_P(\nu_j) =0$. Then    
we obtain $q = (0,0)$, $\overline{p_0} = (1,0)$, and 
the $2$-dimensional lattice polygon 
$P_\sigma = \pi_\sigma(P)$ is contained in the angle 
\[ A_\sigma = \{ (x_1, x_2) \in \R^2 \, \mid \, x_1 \geq 0, x_1 + k x_2  \geq 0 \} \]
and $P_\sigma$ has lattice vertices 
on both sides of $A_\sigma$. Moreover, $P_\sigma$ contains the lattice 
point $(1,0) = \overline{p_0}  = \pi_\sigma(p_0) \in 
\pi_\sigma(F(P))$ in its interior. Elementary convexity considerations show that all above conditions for $P_\sigma$ together with the condition $(0,0) \not\in P_\sigma$ can be satisfied  only if $k =1$ and $(0,1), (1,-1)$ are vertices
of $P_\sigma$. In the last situation the 
lattice vector $2e_1+ e_2 = e_1 + (e_1 + e_2) = \nu_i + \nu_j$ must belongs to  $S_F(P)$, because the line $1 = 2x_1 + x_2$ defines 
a side of $P_\sigma$ and the line $2 = 2x_1 + x_2$ 
contains $(1,0) = \pi_\sigma(p_0) \in \pi_\sigma(F(P))$. The conclusion that   $n_i, n_j, \nu_i + \nu_j \in S_F(P)$ leads to contradiction, since $S_F(P) = S_F(C(P))$, $F(P) = F(C(P))$ and two functions $\ord_{F(P)}$, $\ord_{C(P)}$ are linear on the cone $\sigma \in \widetilde{\Sigma}$  spanned by $n_i$ and $n_j$  (see Remark \ref{prop-supportFP}). So we obtain  
that $(0,0)$ must be a vertex of $P_\sigma$, and   
there exists a vertex $v \in P \cap L_\sigma^0$ such that $\pi_\sigma(v) = (0,0)$. 
 \hfill  $\Box$
\medskip

\noindent
{\bf Proof of Theorem \ref{Mink-sum}.}  Let $P$ be a full dimensional canonically closed Delzant polytope. By \ref{can-closed>1}, we have $C(\lambda P) = \lambda P$ $\forall \lambda \geq 1$. In the proof of \ref{can-sing1} we have shown that 
$F(\lambda P) = F(P) + (\lambda -1)P$ for sufficiently small 
$\varepsilon = \lambda -1 >0$ and the normal fan $\Sigma$ of 
$F((1 + \varepsilon) P)$ equals the normal fan $\widetilde{\Sigma}$ of $F(P) + P$, i.e., $\widetilde{\Sigma}$ is a refinement 
of the normal fan $\Sigma_P$ such that the  corresponding 
refinement of any full dimensional cone $\sigma  \in \Sigma_P(d)$
is determined by the full dimensional domains 
of linearity in $\sigma$ of the $\widetilde{\Sigma}$-piecewise linear 
function $$\delta_P(y) = \ord_{F(P)}(y) - \ord_P(y) = 
\ord_{F(P)}(y) - \< p_\sigma, y \>, \;\; \forall y \in \sigma,$$ 
where $p_\sigma$ is a vertex of $P$ corresponding to 
$\sigma \in \Sigma_P(d)$. 
Since $F(P) \subset p_\sigma + \Theta_\sigma^*$ (see \ref{FP-refinement}), we obtain that $\gamma_\sigma(y) \leq \ord_{F(P)}(y) - \< p_\sigma, y \>$ $\forall y \in \sigma$.  This implies
\[ \Theta_\sigma = \{ y \in \sigma \mid \gamma_\sigma(y) \geq 1 \} 
 \subseteq 
 \{ y \in \sigma \mid \ord_{F(P)}(y) - \< p_\sigma, y \>  \geq 1 \}. \]
 Since $\widetilde{\Sigma}[1] \subset S_F(P)$ and $\ord_{F(P)}(\nu) - \< p_\sigma, \nu \> = 1$ for all $\nu \in S_F(P)$, we obtain 
 the opposite incusion 
\[   
\{ y \in \sigma \mid \ord_{F(P)}(y) - \< p_\sigma, y \>  \geq 1 \} \subset \Theta_\sigma. \]
Therefore, we obtain the equality of the convex sets
\[   
\{ y \in \sigma \mid \ord_{F(P)}(y) - \< p_\sigma, y \>  \geq 1 \} =  \Theta_\sigma \]
for any $\sigma \in \Sigma_P(d)$, i.e., the refinement 
of $\Sigma_P$ defined by the fan $\widetilde{\Sigma}$ is canonical. The equality $\{ y \in \sigma \mid \ord_{F(P)}(y) - \< p_\sigma, y \>  \geq 1 \} =  \Theta_\sigma$ shows that full 
dimensional subcones in the canonical refinement of a cone 
$\sigma \in \Sigma_P$ $1$-to-$1$ correspond to 
rational vectors $\mu  = f - p_{\sigma}  \in {\mathcal M}(\sigma)$, 
where $f $ is a vertex of $F(P)$. 

It remains to show $F(\lambda P) = F(P) + (\lambda -1)P$
for all $\lambda >1$. By \ref{gen-fp}, the $1$-dimensional 
primitive lattice generators of the normal fan of $F(\lambda P)$ 
are contained in $\Sigma_{\rm can}[1] = \{ \nu_1, \ldots, \nu_k\}$, where $\Sigma_{\rm can} = \widetilde{\Sigma}$ is the canonical refinement of the normal 
fan of $\lambda P$. Therefore, we can obtain $F(\lambda P)$ as 
\[ F(\lambda P) = \{ x \in M_\R \mid  \< x , \nu \> \geq \ord_{\lambda P}(\nu ) + 1 \; \; \forall \nu \in 
\widetilde{\Sigma}[1]\}. \]  
On the other hand, the Minkowski sum $F(P) + (\lambda -1)P$ 
corresponds to the sum of two nef $\R$-Cartier 
divisors on $\widetilde{V}$:
 \[  - \sum_{i =1}^k (1 +\ord_P(\nu_i))D_i, \; 
 \; - \sum_{i =1}^k (\lambda -1) \ord_P(\nu_i))D_i. \]
By adding this divisors, we obtain the ample $\R$-Cartier divisor 
on $\widetilde{V}$:
\[  - \sum_{i =1}^k (1 + \lambda \ord_P(\nu_i))D_i = - \sum_{i =1}^k (1 +  \ord_{\lambda P}(\nu_i))D_i  . \]
This implies $F(\lambda P) = F(P) + (\lambda -1)P$. 
\hfill $\Box$ 

\begin{coro}
Let $P$ be a full dimensional generalized Delzant polytope 
$(F(P) \neq \emptyset)$. Then  
\[ \widetilde{P} :=  F(P) + C(P) = F(2C(P)). \]
\end{coro}

\noindent
{\bf Proof.} It is sufficient to apply Theorem \ref{Mink-sum}  to 
the canonically closed polytope $C(P)$ and put $\lambda =2$, because 
$F(C(P)) = F(P)$ (see \ref{can-ord}(b)). 
\hfill $\Box$ 

\begin{coro} \cite[Cor.2.4.3, p.52]{Koe91}\footnote{This result 
of Robert Koelman  was  pointed out to me by Dimitrios Dais.} \label{2-dim-can}
Let $P$ be a $2$-dimensional lattice polytope. Assume that 
 $I(P):= {\rm Conv}(P^\circ \cap M)$
is nonempty. Then the normal fan of the Minkowski sum 
$P + I(P)$ defines a projective toric surface $\widetilde{V}$
with at worst Gorenstein quotient $A_n$-singularities. 
\end{coro}

\noindent
{\bf Proof.} By \ref{int-latt} and \ref{cp-d=2}, we have  $F(P) = I(P)$
and $C(P) = P$. Therefore, 
\[ \widetilde{P} = C(P) + F(P) = P + I(P). \]
Now the statement follows from \ref{can-sing1} and 
\ref{can2}, because $2$-dimensional canonical 
toric singularities are exactly $A_n$-singularities \cite[Prop.11.2.8]{CLS11}. 
 \hfill $\Box$
\medskip

\begin{center}
\begin{tikzpicture}

\begin{scope}[xshift= -2cm, yshift= -2cm] 
\draw[step=1cm,gray,very thin]
(-3.4,-4.4) grid (4.4,4.4) ;

\draw[color= black, solid, very thick] (0,0) -- (1,0) ;
 \draw [color = black, solid, very thick] (0,0) -- (0,-1) ;
  \draw [color = black, solid, very thick] (0,0) -- (-1,-3) ;
 \draw [color = black, solid, very thick] (0,0) -- (4,1) ;
 \draw [color = black, solid, very thick] (0,0) -- (-1,0) ;
 \draw [color = black, solid, very thick] (0,0) -- (-3,2) ;
  \draw [color = black, solid, very thick] (0,0) -- (-1,1) ; 
 
\draw [dashed] (0,0) -- (-2,1) ;
\draw [dashed] (0,0) -- (-1,-1) ;
\draw [dashed] (0,0) -- (-1,-2) ;
\draw [dashed] (0,0) -- (1,1) ;
\draw [dashed] (0,0) -- (2,1) ;
\draw [dashed] (0,0) -- (3,1) ;
\draw [dashed] (0,0) -- (0,1) ;

\draw (-3,2) -- (-1,1) -- (4,1) -- (1,0) --
(0,-1) --(-1,-3) -- (-1,0) -- (-3,2) ;

\fill [color =black] (0,-1) circle (3pt);
\fill [color = black] (1,0) circle (3pt); 
\fill [color = black] (-1,-3) circle (3pt);
\fill [color = black] (-1,0) circle (3pt);
\fill [color = black] (-1,1) circle (3pt);
\fill [color = black] (4,1) circle (3pt);
\fill [color = black] (-3,2) circle (3pt);

\fill (0,1) circle (2pt);
\fill (1,1) circle (2pt); 
\fill (2,1) circle (2pt);
\fill (3,1) circle (2pt);
\fill (-1,-1) circle (2pt);
\fill (-1,-2) circle (2pt);
\fill (-2,1) circle (2pt);

\end{scope} 

\begin{scope}[xshift= -10cm] 
\draw[step=1cm,gray,very thin]
(-2.4,-6.4) grid (4.4,2.4) ;

\draw[ solid, very thick] (-2,2) -- (4,0) -- (0,-6) -- (-2,2); 

\draw[ fill=gray!10!white, very thin] (-1,1) -- (0,1) -- (3,0) 
-- ( 3,-1) -- (1,-4) -- (0,-5) -- (-1,-1) -- (-1,1) ;

 \draw[fill= gray!35!white, solid, very thick] (2,0) -- (3,0) --
 (3,-1) -- (2,-2) -- (2,0); 
 
 \draw[fill= gray!35!white, solid, very thick] (0,-3) -- (1,-3) --
 (1,-4) -- (0,-5) -- (0,-3);

\draw[solid, very thick] (0,1) -- (3,0) ;
\draw[solid, very thick] (3,-1) -- (1,-4); 

\draw[solid, very thick] (0, -5) -- (-1,-1) ;

\draw[fill=gray!35!white, solid, very thick] (-1,1) -- (0,1) --
 (0,0) -- (-1,-1) -- (-1,1);

\draw[fill=gray!35!, solid, very thick]
(-1,1) -- (2,0) -- (0,-3) -- (-1,1);

 \draw[fill=gray!70!white, solid, very thick] (0,0) -- (1,0) --
 (1,-1) -- (0,-2) -- (0,0);

\fill (0,0) circle (3pt);
\fill (1,0) circle (2pt);
\fill (2,0) circle (3pt);
\fill (-1,1) circle (3pt);
\fill (0,-2) circle (2pt);
\fill (0,-3) circle (3pt);
\fill (1,-1) circle (2pt);

\fill (0,1) circle (3pt);
\fill (-1,-1) circle (3pt);
\fill (4,0) circle (3pt);
\fill (-2,2) circle (3pt);
\fill (0,-6) circle (3pt);

\fill (0,-5) circle (3pt);
\fill (1,-4) circle (3pt);

\fill (3,-1) circle (3pt);
\fill (2,-2) circle (2pt);
\fill (1,-3) circle (2pt);
\fill (3,0) circle (3pt);

\draw [ -> ] (-2, 2) -- (-1.1, 1.1) ;
\draw [ -> ] (-2, 2) -- (-0.1, 1.1) ;
\draw [ -> ] (-2, 2) -- (-1.1, -0.9) ;

\draw [ -> ] (4, 0) -- (3.2, 0) ;
\draw [ -> ] (4, 0) -- (3.1, -0.9) ;

\draw [ -> ] (0, -6) -- (0, -5.2) ;

\draw [ -> ] (0, -6) -- (0.9, -4.2) ;

\node [right] at (0,-0.2) {{0}}; 

\draw [dashed] (0,0) -- (0,2/3); 
\draw [dashed] (0,0) -- (-3/5,-3/5); 

\end{scope}
\end{tikzpicture}
\end{center}

\section{Canonical models of hypersurfaces and the  adjunction}

Let $V = V_\Sigma$ be any $d$-dimensional 
normal toric variety defined by a fan $\Sigma \subset N_\R$. 
We consider the Zariski closure $\overline{Z}$ in $V$ of 
an affine hypersurface $Z \subset \T^d$ defined 
by zeros of a Laurent polynomial $f$. 
We need the following  general statement:

\begin{prop} \label{div-hyper}
Let $Z \subset \T^d$ be an 
arbitrary (not necessarily nondegenerate) 
affine hypersurface defined 
by a Laurent polynomial $f({\bf t})$ 
with the Newton polytope 
$P$.  Then  the Zariski closure $\overline{Z}$  of $Z$ in a  
normal toric variety $V_\Sigma$ 
 is linearly equivalent to the torus invariant 
Weyl divisor
\[ D_f:=  - \sum_{\nu \in \Sigma[1]}  \ord_P(\nu) D_\nu.  \]
\end{prop}
 
\noindent
{\bf Proof. } The Weyl divisor 
class group of a normal 
toric variety $V_\Sigma$ does not depend 
on torus orbits of codimension $\geq 2$. Therefore, 
we can assume that the fan $\Sigma$ consists only of cones 
of dimension $0$ and $1$, where  $1$-dimensional cones
$\sigma = \R_{\geq 0} \nu \in \Sigma[1]$ are spanned by 
primitive lattice vectors $\nu \in  \Sigma[1]$.  
The toric variety 
$V_\Sigma$ is smooth, because $V_\Sigma$ is normal.  Therefore, 
Weyl divisors $D$ on $V_\Sigma$ can be identified with  the 
invertible sheaves  ${\mathcal O}(D)$  using the local equations
of $D$ in  $V_\Sigma$. We may consider the Zariski closure $\overline{Z}$ as  
the zero set of a natural global section $s \in H^0(V, 
{\mathcal O}({\overline{Z}}))$ and  use the open affine covering 
\[ V_\Sigma = \bigcup_{ \sigma \in \Sigma(1) } U_\sigma \]
where every open subset $U_\sigma$ consists 
of two torus orbits: the open torus $\T^d$ and the 
closed torus orbit $D_\nu \subset U_\sigma $  isomorphic to  
the divisor $$\{ t_1 = 0 \} \subset 
\C \times (\C^*)^{d-1} = {\rm Spec}\, \C[t_1, t_2^{\pm}, \ldots, t_d^{\pm} ] \cong U_\sigma.$$ 
The maximal ideal $(t_\nu)$ in the local ring of 
$D_\nu \subset U_\sigma \subseteq V_\Sigma$ defines 
a discrete valuation of $\C(t_1, t_2, \ldots, t_d)$ which has value $1$ on $t_\nu := t_1$ and 
value $\ord_P(\nu)$ on the Laurent polynomial $f({\bf t})$. 
 In particu\-lar, the Laurent polynomial 
 $t_1^{-\ord_P(\nu)} f({\bf t})$ defines the 
 local equation of the Zariski closure $\overline{Z}$ in $U_\sigma$. 
 Therefore the rational function 
 $t_\nu^{\ord_P(\nu)}/f({\bf t})$ is the generator of the 
 invertible sheaf ${\mathcal O}(\overline{Z})$ over $U_\sigma$. 
 On the other hand, the invertible sheaf ${\mathcal O}(D_f)$ corresponding to the torus invariant divisor $D_f:=   - \sum_{\nu \in \Sigma[1]}  \ord_P(\nu) D_\nu$ is generated over $U_\sigma$ 
 by $t_\nu^{\ord_P(\nu)}$, where $t_\nu$ is the generator of 
 the vanishing ideal of the closed torus orbit $D_\nu \subset U_\sigma$. The multiplication by the Laurent polynomial $f$ defines 
 an isomorphism ${\mathcal O}(\overline{Z}) \cong {\mathcal O}(D_f)$. Therefore, $D_f$ is linearly equivalent to $\overline{Z}$.
 \hfill $\Box$

\begin{coro}
The divisor class  
$$[ \overline{Z}] \in \bigoplus_{\nu \in \Sigma[1]} \Z D_\nu \cong 
\Z^{\Sigma[1]} $$ 
of the Zariski closure  $\overline{Z}$ of $Z$ in $V_\Sigma$ depends only on the Newton polytope $P$ of the defining  
Laurent polynomial of $f$. If one chooses  
another Laurent polynomial $f({\bf t}) {\bf t}^m$  $(m \in M)$ defining the same affine 
hypersurface  $Z \subset \T^d$, then $P$ will be 
replaced by  the shifted  Newton polytope $P +m$ and 
one obtains 
\[ [ \overline{Z}_{f({\bf t}){\bf t}^m} ] = \sum_{\nu \in \Sigma[1]}  (\ord_P(\nu) + \langle m, \nu \rangle) D_\nu.  \]
\end{coro}

\begin{defi} Let $Z \subset \T^d$ be a nondegenerate affine 
hypersurface defined by a Laurent polynomial $f({\bf t})$ 
with the Newton polytope $P$. We assume that $F(P) \neq \emptyset$ and call the Zariski closure $\widetilde{Z} \subset 
\widetilde{V}$ in the projective toric variety associated 
with the rational polytope $\widetilde{P}:= F(P) + C(P)$ the {\em canonical model} of $Z$. 
\end{defi}

\begin{prop}
The canonical model $\widetilde{Z}$ is a semiample big 
$\Q$-Cartier divisor on toric variety $\widetilde{V}$. 
\end{prop}

\noindent
{\bf Proof.} Since the full-dimensional polytope 
$C(P)$ is a Minkowski summand
of $\widetilde{P}$, we obtain a natural birational 
toric morphism 
\[  \varrho\,:\,  \widetilde{V} \to V_{C(P)}  \]
and the canonical model $\widetilde{Z}$ is the pull-back 
under $\varrho$ of the ample $\Q$-Cartier divisor 
$\overline{Z}_{C(P)} \subset V_{C(P)}$. This implies 
the statement. 
\hfill $\Box$ 
\medskip

\begin{rem}
Note that if $P =C(P)$, then the toric 
morphism  $\varrho\,:\,  \widetilde{V} \to V_{C(P)} =V_P$
is  determined by the canonical refinement $\Sigma^{\rm can}_P$ of the normal 
fan $\Sigma_P$ and the singularities of $\widetilde{Z} $ 
are toroidal. In general, we have $P \neq C(P)$ if 
$P$ contains facets $Q$ whose normals do not belong to $S_F(P)$. These "bad facets" of $P$ define divisors on $V_P$ and 
on  $\overline{Z}$ that must disappear on $\widetilde{V}$ and on the canonical (resp. minimal)
models $\widetilde{Z}$ (resp. $\widehat{Z}$). Therefore, in general, there is no any birational toric morphism 
$$\rho' \; :\; \widetilde{V} \to V_P.$$ This is the reason why
the singularities of canonical and minimal models $\widetilde{Z}$ and $\widehat{Z}$ are not toroidal in general. We illustrate this fact  in Example \ref{simp-minimal} (Case 1).
\end{rem} 

The next theorem is very important. 

\begin{thm} \label{norm-adj}
The canonical model  $\widetilde{Z} $ is a normal variety   and it satisfies the adjunction formula:
\begin{align} \label{adj}
K_{\widetilde{Z}} = (K_{\widetilde{V}} + \widetilde{Z})\mid_{\widetilde{Z}}.
\end{align}
\end{thm}

The main technical difficulty  in the proof of Theorem \ref{norm-adj} is the fact that the hypersurface $\widetilde{Z} \subset
\widetilde{V}$ may contain some torus orbits 
in the toric variety $\widetilde{V}$. This could make problems 
for the adjunction as it can be seen from the next 
example. 

\begin{ex}
Consider a surface $S \subset \P(1,1,a,a) =V$ $(a \geq 1)$ 
defined by weighted homogeneous equation of degree $a+1$: 
\[ z_0 z_2 - z_1 z_3 =0. \] Then the surface
$S$ is smooth and  it is isomorphic to the ruled surface $\F_{a-1}$. 
However, if $a \geq3$ the adjunction formula for $S$ does not hold, because
$\F_{a-1}$ is not a del Pezzo surface.  The adjunction in this case
fails, because  
the singularities of the weighted projective space $V$ include a $1$-dimensional 
torus orbit $\Theta \subset V$ contained in $S$. This is an example of a quasi-smooth divisor $S$ in the weighted projective space that is not well-formed in sense of Fletcher \cite{Fle00}. 
\end{ex}

Our proof of Theorem \ref{norm-adj} is based on the following 
statement:

\begin{prop} \label{codim2-orbits}
None of the $(d-2)$-dimensional torus orbits $T_\sigma \subset \widetilde{V}$ corresponding 
to $2$-dimensional cones $\sigma \in \widetilde{\Sigma}$ is contained in the canonical model $\widetilde{Z}$. 
\end{prop}

\noindent 
{\bf Proof.} Let $\nu_i, \nu_j \in N$ be primitive lattice generators
such that $\sigma = \R_{\geq 0} \nu_i + \R_{\geq 0}\nu_j$. Consider  two supporting hyperplanes 
\[ L_i\; : \;  \langle x, \nu_i \rangle = {\ord}_P(\nu_i) \;\; \text{and} 
\;\;  L_j\; : \;   \langle x, \nu_j \rangle = {\ord}_P(\nu_j). \]
Then $L_i \cap P \neq \emptyset$ and  
$L_j \cap P \neq \emptyset$. Moreover, the Zariski closure $Z_\sigma$ of $Z$ 
in the open affine toric subvariety $U_\sigma \subset \widetilde{V}$
corresponding to the $2$-dimensional cone $\sigma$ 
contains the unique closed torus 
orbit $T_\sigma \subset U_\sigma$ if 
and only if $L_i \cap L_j \cap P \neq \emptyset$. It remains to 
apply \ref{codim2faces}(b).  \hfill $\Box$

\noindent
{\bf Proof of \ref{norm-adj}.} Since a normal toric variety 
$\widetilde{V}$ is Cohen-Macaulay and $\widetilde{Z}$ is 
a hypersurface in $\widetilde{V}$, by Serre's criterion for
normality, 
it is sufficient
to show that $\widetilde{Z}$ is smooth in codimension $1$. The latter enough to check on the open subset  
\[ \widetilde{U}^{(2)} := \bigcup_{\sigma \in \widetilde{\Sigma} \atop \dim \sigma =2} U_\sigma, \]
because the complement $W :=\widetilde{V} \setminus \widetilde{U}^{(2)}$ has codimension $3$ in $\widetilde{V}$. 
By \ref{codim2-orbits}, for every $2$-dimensional cone $\sigma \in \widetilde{\Sigma}$ 
the closed torus orbit $T_\sigma$ is not contained in the Zariski
closure $\widetilde{Z}$. Let $Q:= L_i \cap L_j \cap P$. Then 
$Q$ is a face of $P$. If $\dim Q =0$, then $T_\sigma \cap \widetilde{Z} = \empty$ and therefore 
$U_\sigma \cap \widetilde{Z}$ is
smooth. If $\dim Q > 0$, then it follows  from the nondegeneracy
of $Z$ that the intersection $T_\sigma \cap \widetilde{Z}$ 
is transversal along a smooth reduced divisor in $T_\sigma$, 
and singularities of $U_\sigma \cap \widetilde{Z}$ must 
be contained in  $T_\sigma \cap \widetilde{Z}$, i.e. they
have at least codimension $2$. 

Finally, the adjunction holds for the canonical model $\widetilde{Z} \subset \widetilde{V}$, because $\widetilde{Z}$ is 
transversal to all torus orbit in $U^{(2)}$ and the adjunction holds on the open subset $U^{(2)} \subset \widetilde{Z}$.    
\hfill $\Box$

\section{Minimal models}

\begin{thm} \label{can-models}
The projective hypersurface $\widetilde{Z} \subset \widetilde{V}$ 
is a $\Q$-Gorenstein variety with nef  $\Q$-Cartier canonical
divisor $K_{\widetilde{Z}}$ and with at worst canonical singularities.
\end{thm}

\noindent
{\bf Proof.} 
We take a common regular refinement $\Sigma$ of 
 the normal fan $\Sigma_P$  and of the fan $\widehat{\Sigma}$.
Denote by $\overline{Z}$ the Zariski closure of $Z$ in the 
smooth toric variety $V:= V_\Sigma$. 
By \ref{div-hyper}, we have   $\overline{Z}  = \sum_{\nu \in \Sigma[1]} - 
  \ord_P(\nu) D_\nu$. Together with 
  $K_{{V}} = \sum_{\nu \in {\Sigma}[1]} - D_\nu$  we obtain
\[ K_V + \overline{Z}  = \sum_{\nu \in \Sigma[1]} (- 1 - 
\ord_P(\nu) ) \cdot D_\nu. \]
We may consider $\Sigma$ as a regular refinement of $\widetilde{\Sigma}$ and obtain the birational toric morphism $\varphi\, : \, V \to  \widetilde{V}$.  Then 
\begin{align} \label{disc1}
 K_{{V}} + \overline{Z} =  \varphi^*(K_{\widetilde{V}} + \widetilde{Z})
 + \sum_{\nu \in \Sigma[1] \setminus \widetilde{\Sigma}[1]}  a(\nu) D_\nu. 
 \end{align}
 In order to compute the discrepancies $a(\nu)$ we write the right hand side in the last equality as a linear combination of the toric 
 divisors $D_{\nu}$ ($\nu \in \Sigma[1]$) using 
\begin{align*}
\varphi^*(K_{\widetilde{V}} + \widetilde{Z}) = & \sum_{ \nu \in \Sigma[1]} -\ord_{F(P)}(\nu) D_\nu \\ = & \sum_{\nu \in \widetilde{\Sigma}[1]} - 
\ord_{F(P)}(\nu)D_\nu +  \sum_{ \nu \in \Sigma[1] \setminus \widetilde{\Sigma}[1]}  -\ord_{F(P)}(\nu) D_\nu  \\
= &   \sum_{\nu \in \widetilde{\Sigma}[1]} (-1 - 
\ord_P(\nu))D_\nu +  \sum_{ \nu \in \Sigma[1] \setminus \widetilde{\Sigma}[1]}  -\ord_{F(P)}(\nu) D_\nu, 
\end{align*}  
\[  \varphi^*(K_{\widetilde{V}} + \widetilde{Z})
 + \sum_{\nu \in \Sigma[1] \setminus \widetilde{\Sigma}[1]}  a(\nu) D_\nu =  \sum_{\nu \in \widetilde{\Sigma}[1]} (-1 - 
\ord_P)D_\nu + \sum_{\nu \in \Sigma[1] \setminus 
\widetilde{\Sigma}[1]} ( -  \ord_{F(P)} + a(\nu))  \cdot D_\nu. \]
This implies 
\[ -1 - 
\ord_P(\nu) =   -  \ord_{F(P)} + 
 a(\nu)  
 \;\; \forall \nu \in \Sigma[1] \setminus \widetilde{\Sigma}[1]. \]
Now we restrict equation (\ref{disc1}) to $\overline{Z}$ and 
obtain 
\begin{align} \label{disc2}
 K_{\overline{Z}} =  \varphi^*(K_{\widetilde{Z}})
 + \sum_{\nu \in \Sigma[1] \setminus \widetilde{\Sigma}[1]}  a(\nu) (D_\nu \cap \overline{Z}). 
 \end{align}
Since  
\begin{align} \label{disc3} 
 a(\nu) =  \ord_{F(P)}(\nu)  -  \ord_P(\nu) -1 \geq 0 \; \; 
\forall \nu \in N \setminus \{ 0 \},  
 \end{align}
the singularities of $\widetilde{Z}$ are at worst canonical.  

\hfill $\Box$

\begin{thm} \label{min-models}
Let $\widehat{\Sigma}$ be a maximal 
projective simplicial refinement 
of the normal fan $\widetilde{\Sigma} = \Sigma_{\widetilde{P}}$ with 
$\widehat{\Sigma}[1] = S_F(P)$. Then the  Zariski closure $\widehat{Z}$ of the affine hypersurface 
$Z \subset \T^d$ in  $\widehat{V}:= V( \widehat{\Sigma})$ is 
a projective minimal model of  $Z$, i.e., $\widehat{\Sigma}$
is a $\Q$-factorial projective algebraic variety with at worst 
terminal singularities and with semiample canonical class 
$K_{\widehat{Z}}$. 
Moreover, the crepant birational toric morphism $\varphi\, :  \widehat{V} \to \widetilde{V}$ 
induces a 
birational crepant 
morphism $\varphi\,: \, \widehat{Z} \to \widetilde{Z}$. 
\end{thm}

\noindent
{\bf Proof.} 
Therefore, the same arguments as above show that the singularities of $\widehat{Z}$ are at worst terminal and one the induced birational morphism 
\[ \varphi\,: \, \widehat{Z} \to \widetilde{Z} \]
is crepant, because  the equality 
\[ a(\nu) =  \ord_{F(P)}(\nu)  -  \ord_P(\nu) -1 =0 \]
holds if and only if $\nu \in S_F(P)$. 
\hfill  $\Box$ 

Let us consider some simplest examples in arbtrary dimension $d$ showing the roles of $F(P)$, $S_F(P)$ and 
$C(P)$ in constructing minimal models of non-degenerate toric hypersurfaces $Z \subset \T^d$. 

\begin{ex}\label{simp-minimal}
Let $P \subset \R^d$ be the $d$-dimensional lattice polytope defined as 
\[ P := \{ x \in \R^d_{\geq 0} \; \mid \; a \leq \sum_{i =1}^d x_i \leq b \}  \]
for some integers $0 \leq a < b$. Then 
\[ F(P)=  \{ x \in \R^d_{\geq 1} \; \mid \; a +1  \leq \sum_{i =1}^d x_i \leq b -1 \}. \] 
So  $F(P) \neq \emptyset$ if and only if $b \geq \max \{d+1, a+2\} $. 

Let $e_1, \ldots, e_d$ be the standard basis of $\Z^d \subset \R^d$. There are two possibilities for $S_F(P)$ if $F(P) \neq \emptyset$:  

{\bf Case 1.} $S_F(P) = \{ e_1,  \ldots, e_d, - \sum_{i=1}^d e_i \}$. This happens if and only if $b \geq d+1$ and $a +2 \leq d$. 
In this case, we have 
\[ C(P) =\{ x \in \R^d_{\geq 0} \; \mid \;  \sum_{i =1}^d x_i \leq b \}, \]
and the minimal model $\widehat{Z} = \widetilde{Z}$ is a  projective hypersurface of degree $b$  
in $\P^d$ that may have an isolated terminal singularity at the  origin $0 \in \C^d \subset \P^d$ if $a \geq 2$. In the latter case,   $C(P) \neq P$ and
the isolated singularity of the minimal model  is not toroidal as soon as $a \geq 3$.  

{\bf Case 2.} $S_F(P) = \{ e_1,  \ldots, e_d, - \sum_{i=1}^d e_i ,  \sum_{i=1}^d e_i  \}$. This happens if and only if  $b \geq a +2 \geq d+1$. 
In this case, we have $C(P) = P$ and  
the minimal model  
$\widehat{Z} = \widetilde{Z} = \overline{Z}$ is a smooth projective hypersurface, the  Zariski closure of $Z$ in the smooth toric variety $\widehat{V} = \widetilde{V} = V_P$ obtained from $\P^d$ by the blow-up of the origin 
$0 \in \C^d \subset \P^d$. 
 
\end{ex}

\section{The Kodaira dimension of toric hypersurfaces}

We begin this section with simple illustrating 
examples showing some 
properties of the Kodaira dimension and of the Iitaka fibration 
for $(d-1)$-dimensional smooth projective 
hypersurfaces $X_{a,b}$ of bi-degree $(a,b)$ in the product 
$\P^{d_1} \times \P^{d_2}$ of two projective spaces  $(d= d_1 + d_2)$. We use these 
examples to illustrate general properties of the Kodaira 
dimension and the Iitaka fibration 
for canonical models of nondegenerate affine 
hypersurfaces $Z \subset \T^d$ defined by Laurent 
polynomials with $d$-dimensional Newton 
polytopes $P$. 

The hypersurfaces $X_{a,b} \subset \P^{d_1} \times \P^{d_2}$ have nonnegative Kodaira dimension $\kappa (X_{a,b})$  
if and only if 
$a \geq d_1 +1$ and $b \geq d_2 + 1$.  

\begin{itemize}

\item 
If $a> d_1 + 1$ and 
$b > d_2 +1$, then the canonical class of $X_{a,b}$ is ample 
and $\kappa(X_{a,b}) = \dim X_{a,b} = d-1$, i.e., 
$X_{a,b}$ is a smooth projective 
variety of general type. 

\item
If $(a, b)= 
(d_1 + 1, d_2 +1)$,  then $X_{a,b}$ has trivial  
canonical class and $X_{d_1+1,d_2+1}$ is a 
$(d-1)$-dimensional  Calabi-Yau hypersurface, i.e., $\kappa(X_{d_1+1, d_2+1}) = 0$. 

\end{itemize}

A more interesting situation appears if e.g.  
$a = d_1 +1$, but $b > d_2 +1$. In this case the canonical 
invertible sheaf on $X_{a,b}$ is the restriction to  $X_{a,b}$  of the semiample sheaf  
${\mathcal K} := {\mathcal O}(0, b - d_2 -1)$ on 
$\P^{d_1} \times \P^{d_2}$. The global sections of ${\mathcal K}$ 
define the natural surjective 
projective morphism $\vartheta\, :\, \P^{d_1} \times 
\P^{d_2} \to \P^{d_2}$ which is a trivial $\P^{d_1}$-fibration,  
a  Fano-fibration, over $\P^{d_2}$. The restriction of $\vartheta$
to the hypersurface $X_{a,b}$ defines a proper surjective 
morphism $X_{a,b} \to \P^{d_2}$ whose general fibers 
are anticanonical hypersurfaces of degree $d_1 +1$ in 
$\P^{d_1}$. 

There are two cases:

\begin{itemize}

\item 
If $d_1 =1$, then the anticanonical hypersurface in $\P^{d_1}$ 
consists of two points in $\P^1$ and the morphism 
$\vartheta\, :\, X_{2,b} \to \P^{d_2}$ is a double covering 
of $\P^{d_2}$, 
a higher dimensional analog of hyperelliptic curves. In 
particular, the Kodaira dimension of $X_{2,b}$ is 
$d-1$.  

\item
If $d_1 > 1$, then for a general point $p \in \P^{d_2}$ 
the  fiber $\vartheta^{-1}(p) \cap X_{d_1 +1,b}$ is an irreducible 
$(d_1 -1)$-dimensional Calabi-Yau 
hypersurfaces of degree $d_1 + 1$ 
in $\P^{d_1}$.  By \ref{product}, the Newton polytope $P$ 
of $X_{d_1 +1, b}$ is a product of two simplices 
of dimensions $d_1$ and $d_2$.   The Fine 
interior $F(P)$ is a $d_2$-dimensional simplex.  It is easy to see that the pluricanonical ring of $X_{d_1+1, b}$ is the $(b-d_2-1)$-Veronese ring of $\P^{d_2}$ 
and that 
the Kodaira dimension of $X_{d_1 +1, b}$ equals 
$d_2$, i.e., 
\[ \kappa(X_{d_1 +1, b}) = \dim F(P) = d_2 < d-1. \]  
\end{itemize}
\medskip

Now we consider a general case of the  
toric fibration $\vartheta\,:\,  \widetilde{V} \to 
V_{F(P)}$ corresponding to the Minkowski summand $F(P)$ of the 
full dimensional polytope $\widetilde{P} = F(P) + C(P)$, where 
$k:= \dim F(P) < d$.  Define the  $(d-k)$-dimensional  sublattice 
 \[ N^{F} := \{ n \in N \mid  \< x, n \> = \ord_{F(P)}(n) 
 \;\; \forall  x \in F(P) \}  \subset N. \]
 The sublattice  $N^F \subset N$ is a direct summand of $N$ 
of rank $d-k$,  i.e., there exists a complementary    
sublattice $N_F$ of rank $k$ such that $N = N_{F} \oplus N^F$ and a natural surjective homomorphism 
\[\pi_F\, : \,  N \to N/N^F \cong N_F.    \]
Denote by $M_F \subset M$ the $k$-dimensional 
orthogonal complement 
of  $N^F$ in $M$. 
We obtain  the natural surjective homomorphism
\[ \pi^F\, :\, M \to M^F:= M/M_F  \]
and denote by $P^F$   the $(d-k)$-dimensional lattice  polytope 
$\pi^F(P) \subset M^F_\R$, i.e.,  the $\pi^F$-projection 
 of the lattice polytope $P$ onto $(d-k)$-dimensional lattice polytope 
 $P^F \subset M^F_\R$. 

 \begin{prop}\label{fibration1}
 The $(d-k)$-dimensional 
 lattice polytope $P^F \subset M_\R^F$ has the following properties$:$  
 
 {\rm (a)}  $F(P^F)= \pi^F(F(P))$, i.e., the Fine interior 
 of the lattice polytope $P^F$ is the rational point $\pi^F(F(P)$; 

{\rm (b)} $S_F(P^F) = N^F \cap S_F(P) \subset N^F$; 

{\rm (c)} $\Phi^F:= {\rm Conv} \left(N^F \cap S_F(P) \right)$ is a 
$(d-k)$-dimensional canonical Fano polytope, i.e., a $(d-k)$-dimensional lattice polytope containing only the origin 
$0 \in N^F$ as $N^F$-lattice point in the relative 
interior of  $\Phi^F$. 
\end{prop}  

\noindent
{\bf Proof.} The lattice $M^F$ is dual to $N^F$. Therefore, 
\[  F(P^F) = \{ x \in M_\R^F \; |\; \langle x, \nu \rangle\geq \ord_{P^F}(\nu) + 1 \;\; \forall \nu \in N^F \setminus \{ 0 \} \}. \] 
Since $P^F$ is $\pi^F$-projection of $P$ and any linear function 
$\langle *, \nu \rangle$ with $\nu \in N^F$ is constant 
on fibers of $\pi^F$-projection, we obtain that 
$\ord_P(\nu) = \ord_{P^F}(\nu)$ for all $\nu \in N^F$. 
This implies  that $p^F:= \pi^F(F(P))$ is the Fine interior of $P^F$. \hfill  $\Box$

\begin{thm}
Let $Z \subset \T^d$ be a nondegenerate affine 
toric hypersurface defined by a Laurent polynomial $f$ 
with a  Newton polytope $P$. If $k := \dim F(P) \geq 0$, 
then the Kodaira dimension $\kappa(\widetilde{Z})$ of the canonical model 
$\widetilde{Z}$ equals 
\[ \kappa(\widetilde{Z}) = \min \{ k, d-1 \} \]
and we have the following three cases: 

{\rm (a)} If $k = d$,  then the Iitaka fibration $\widetilde{Z} \to V_{F(P)}$ is birational on its image.   

{\rm (b)}  If $k =d-1$, then  the minimal model
 $\widetilde{Z}$ is birational to a double cover
of a $(d-1)$-dimensional toric variety $V_{F(P)}$, i.e., $\widetilde{Z}$ is a higher dimensional analog of hyperelliptic curves of genus $g \geq 2$.

{\rm (c)} If $0 \leq k < d-1$, then the Iitaka fibration 
\[   \widetilde{Z} \to V_{F(P)}  \]
is induced by  the canonical toric $\Q$-Fano fibration $\vartheta\, :\, \widehat{V} \to V_{F(P)}$ 
whose generic fiber is  
isomorphic to a nondegenerate irreducible 
$(d-1-k)$-dimensional 
hypersurface of Kodaira dimension $0$ 
in some toric $\Q$-Fano variety defined by the $(d- k)$-dimensional lattice polytope $P^F$ with $0$-dimensional Fine interior.  
\end{thm}

\noindent
{\bf Proof. } 
By the adjunction formula, $(K_{\widetilde{V}} + \widetilde{Z})|_{\widetilde{Z}}$ is the canonical class $K_{\widetilde{Z}}$. 
So we obtain the linear maps 
\[ \Psi_m \; : \; H^0(\widetilde{V}, {\mathcal O}( m(K_{\widetilde{V}} + \widetilde{Z})) \to  H^0(\widetilde{Z}, {\mathcal O}( m K_{\widetilde{Z}})), \;\; m \geq 0.   \]
Consider the toric morphism $\vartheta\, :\, \widetilde{V} \to V_{F(P)}$. Then   $K_{\widetilde{V}} + \widetilde{Z} = 
\vartheta^* L$ for some 
ample divisor $L$ on the toric variety $V_{F(P)}$. Therefore, 
the dimensions $h^0 (\widetilde{V}, {\mathcal O}( m(K_{\widetilde{V}} + \widetilde{Z}))$ grow as degree $\dim F(P)$ polynomial of $m$. 

We consider two cases: 

Case 1. $k:= \dim F(P) =  d$. Then the dimensions  $h^0 (\widetilde{V}, {\mathcal O}( m(K_{\widetilde{V}} + \widetilde{Z}))$
grow as degree $d$ polynomial of $m$. This implies that the 
restrictions of these global sections to hypersurface 
$\widetilde{Z}$ grow as at least degree $d-1$ polynomial of $m$. 
The latter implies that only $\kappa(\widetilde{Z}) = \dim \widetilde{Z} = d-1$ is possible. 

Case 2. $k:= \dim F(P) < d$.  We claim that in this case 
all maps $\Psi_m$ $(m \geq 0)$ 
are injective (the latter implies $\kappa(Z) \geq k$).

Indeed, 
any global section   
$$s \in H^0(\widetilde{V}, {\mathcal O}( m(K_{\widetilde{V}} + \widetilde{Z}))$$ is represented by a  
Laurent polynomial $g({\bf t})$ whose Newton polytope 
$P(g)$ 
is contained in the  $k$-dimensional 
polytope $m F(P)$.  If the restriction of such $g$  
to ${Z}$ is zero, then the Laurent polynomial 
$f$   divides $g$ and the Newton 
polytope $P = Newt(f)$ can be embedded into the Newton polytope 
of $Newt(g)$. 
The latter is impossible by dimension reasons unless 
$g = 0$. 

In order to get the opposite inequality $\kappa(\widetilde{Z}) 
\leq k$ we remark
that on the toric variety $\widetilde{V}$ the nef-divisor  $ K_{\widetilde{V}} + \widetilde{Z}$ 
 defines 
a toric morphism 
\[ \varphi \, : \, \widetilde{V} \to V_{F(P)}:=
{\rm Proj}  \bigoplus_{m \geq 0} H^0(\widetilde{V}, {\mathcal O}( m(K_{\widetilde{V}} + \widetilde{Z}))  \]
and  $ V_{F(P)}$ is a toric variety of dimension $k$. 

The fibers of $\varphi$ over the dense torus 
orbit $U \subset  V_{F(P)}$ are $(d-k)$-dimensional 
canonical $\Q$-Fano toric varieties and 
the restiction of global sections of $ {\mathcal O}( m(K_{\widetilde{V}} + \widetilde{Z}))$  are 
trivial. Therefore, the $\varphi$-images 
of the generic intersections of $  \widetilde{Z}$ with  
these fibers are $(d-k-1)$-dimensional. These irreducible 
$(d-k-1)$-dimensional subvarieties of $  \widetilde{Z}$ are 
mapped to points in $V_{F(P)}$. Therefore, the dimension
of the pluricanonical image of $  \widetilde{Z}$ can be 
at most $(d-1) - (d-k-1) =k$. Therefore, we obtain 
$\kappa(   \widetilde{Z}) \leq k$.  
\hfill  $\Box$

Consider an illustrating example:
 
\begin{ex} \label{enri}
Let $P = C(P)$ be canonically closed  
$3$-dimensional lattice simplex with vertices 
\[ (0,0,0), (3,0,0), (1,3,0), (2,0,3) \in M_\R = \R^3. \]
The simplex $P$ contains no lattice points in its interior, but 
the Fine interior $F(P)$ is $1$-dimensional rational segment of length $1/3$: 
\[ F(P) = \left[ (4/3,1,1), (5/3, 1,1)     \right] \subset M_\R.    \]
The $1$-dimensional sublattice $M_F \subset M$ is spanned by 
$(1,0,0)$. The projection $P^F \, : \, M \to M^F = M/M_F \cong \Z^2$ is the map $(m_1, m_2,m_3)  \mapsto (m_2,m_3)$. 
The image  $P^F(P)$ of $P$ under this projection  is the reflexive lattice triangle with vertices  $(0,0), (0,3), (3,0)$  having 
the $0$-dimensional Fine interior $\{ (1,1) \}$. 
One can consider the nondegenerate affine surface $Z \subset \T_3 \cong (\C^*)^3$ defined by the equation 
\[ f({\bf t}) = 1 + t_1^3 + t_1t_2^3 + t_1^2t_3^3 =0 \]
with the Newton polytope $P$. 
The canonical model $\widetilde{Z}$ is an elliptic surface of Kodaira dimension $1$ having the natural surjective morphism $\widetilde{Z} \to \P^1$, if we take   $t_1$ as 
local affine coordinate on $\P^1$. Thus we obtain a  
$1$-parameter $t_1$-family of plane elliptic cubic curves  determined by  the affine $(t_2,t_3)$-equations
\[ (1 + t_1^3) + t_1t_2^3 + t_1^2 t_3^3 = 0,  \;\; t_1 \in \C \]
having  the reflexive Newton polygon with 3 vertices 
$(0,0), (3,0), (3,0)$.

\begin{center}
\begin{tikzpicture}[scale=0.5]

\fill (0,3.3) circle (3pt);
\fill (8.7, 0) circle (3pt);
\fill (5.8, 11) circle (3pt);

\draw[ solid, very thick] (0,3.3) -- (8.7, 0) -- (5.8, 11) -- (0,3.3)  ;

\fill (5.8, 11) circle (3pt);

\fill (8.9, 6.1) circle (3pt);

\draw[dashed, very thick] (0,3.3) --  (8.9, 6.1) ; 

\draw[ solid, very thick]  (8.9, 6.1)  -- (8.7, 0)   ; 

\draw[ solid, very thick]  (8.9, 6.1)  -- (5.8, 11)   ; 

\fill [color = gray] (4.9, 6.8) circle (3pt); 

\fill (7.8, 5.7) circle (3pt); 

\fill [color=black] (5.866, 6.434) circle (4pt); 

\fill [color = black] (6.834, 6.066) circle (4pt); 

\draw[ color= black, solid, very thick] (5.866, 6.434) -- 
(6.834, 6.066) ; 

\draw [very thin, dashed] (6.834, 6.066) -- (7.8, 5.7) ; 
\draw[very thin] (7.8, 5.7) -- (10.7, 4.6);

\draw[very thin, dashed] (5.866, 6.434) -- ( 4.9, 6.8);

\draw[very thin, dashed] (4.9, 6.8) -- (2,  7.9); 

\fill (2, 7.9) circle (3pt); 

\draw[very thin] (2.968, 7.532) --(2, 7.9)  ;

\node [left] at (0,3.3) {{ $(0,0,0)$}}; 
\node [right] at (8.7,0) {{ $(3,0,0)$}}; 
\node [right] at (8.9, 6.1) {{ $(1,3,0)$}}; 
\node [right] at  (5.8,11) { $(2,0,3)$}; 

\draw[very thin] (5.8, 11) -- (5.8, 1.1); 
\draw[very thin, dashed] (8.9, 6.1) -- (2.9, 2.2);

\draw[ - > , very thin] (0, 3.3) -- (0,13.2); 

\fill (0, 6.6) circle (3pt);
\fill (0, 9.9) circle (3pt);
\fill (0, 13.2) circle (3pt);
\fill [color = gray] (6.0, 7.2) circle (3pt);
\draw[very thin, dashed] (8.9,6.1) --  (6.0, 7.2) ;
\fill [color = gray] (4.0, 5.9) circle (3pt);
\fill [color = gray] (2.0, 4.6) circle (3pt);
\draw[very thin] (0,13.2) --  (4.25, 8.95) ;
\draw[very thin, dashed] (6.0,7.2) --  (4.25, 8.95) ;
\draw[very thin] (0,13.2) --  (5.8, 11) ;
\fill  (4.0, 9.2) circle (3pt);
\fill  (2.0, 11.2) circle (3pt);

\draw[ ->,  very thin, dashed] (0, 3.3) -- (8 , 8.5) ;

\draw [ - > , very thin, dashed] (8.7, 0) -- (11.6, -1.1) ;  
\end{tikzpicture}
\end{center}

The fan $\widetilde{\Sigma} = \widehat{\Sigma}$ defining 
$3$-dimensional toric variety $\widetilde{V}$ with at worst canonical singularities is generated by $5$ 
lattice vectors 
\[ S_F(P) = \widetilde{\Sigma}[1] = 
\widehat{\Sigma}[1]  = \{ (0,1,0), (0,0,1), (-1,-1), (3,-1,-2), (-3,-2,-1) \subset N_\R. \]
The terminal toric variety $\widehat{V}$ equals  $\widetilde{V}$. It has $6$ isolated terminal $\mu_3$-quotient  singularities and it is generically 
a toric $\P^2$-fibration over $\C^*$. 
The $2$-dimensional sublattice $N^F = {\rm Span}( (0,1,0), (0,0,1) ) \subset N$ contains altogether three lattice vectors from $S_F(P)$ spanning the fan of $\P^2$:  
\[  \{ (0,1,0), (0,0,1), (0,-1,-1)  \} \subset N^F. \]

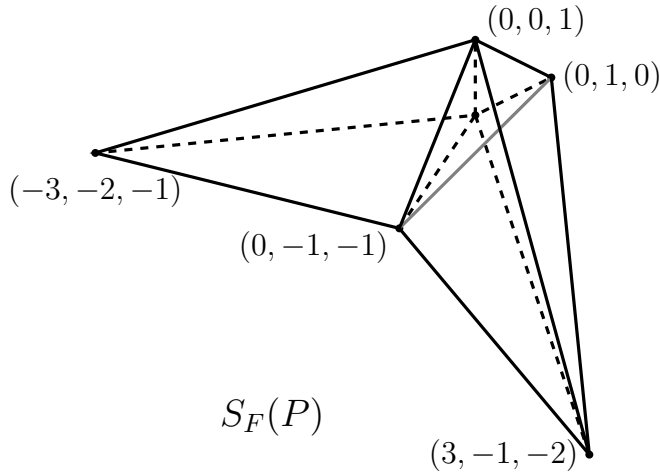
 \begin{center}
\begin{tikzpicture}[scale=0.5]
\fill (4,4) circle (3pt);
\fill (0,0) circle (3pt);
\fill (5,-6) circle (3pt);
\fill (-8,2) circle (3pt);
\fill (2,3) circle (3pt);
\fill (2,5) circle (3pt);

\draw[very thick, dashed] (5, -6) -- (2,3) -- (-8, 2);
\draw[very thick, dashed] (0, 0) -- (2,3) -- (4, 4);
\draw[very thick, dashed]  (2,3) -- (2, 5);
\draw[color=gray, very thick]  (0,0) -- (4, 4);
\draw[very thick] (0, 0) -- (2,5); 
\draw[very thick] (2, 5) -- (5,-6); 
\draw[very thick]  (5,-6) -- (4,4) --(2,5) -- (-8, 2) -- (0, 0) -- (5,-6) ;

\node[left] at (0,-0.5) { {\bf $(0,-1,-1)$}};
\node at (-8,1) { {\bf $(-3,-2,-1)$}};
\node[left] at (5,-6) { {\bf $(3,-1,-2)$}};
\node[right] at (4,4) { {\bf $(0,1,0)$}};
\node[right] at (2,5.5) { {\bf $(0,0,1)$}};

\node[right]  at (-5,-5)  {\large {$S_F(P)$}};

\end{tikzpicture}
\end{center}

\end{ex}

We end  this section with a combinatorial formula 
for the top intersection number $(K_{\widehat{Z}})^{d- 1}$. 

\begin{thm} 
Let $P$ be a $d$-dimensional lattice polytope with 
$F(P) \neq \emptyset$. Denote by $\widehat{Z}$ a minimal model obtained of nondegenerate hypersurface $Z \subset \T^d$ with the Newton polytope $P$ which is obtained from 
the canonical models $\widetilde{Z}$ by 
a crepant morphism $\varphi\, :\, \widehat{Z} \to 
\widetilde{Z}$. Then 
$$(K_{\widehat{Z}})^{d-1} =(K_{\widetilde{Z}})^{d-1} = \begin{cases} 
{\rm Vol}_d(F(P)) + \sum_{ Q \prec F(P) \atop \dim Q =d-1} {\rm Vol}_{d-1}(Q), & \text{ if $\dim F(P) =d$;} \\
2 {\rm Vol}_{d-1}(F(P)), & \text{ if $\dim F(P) =d-1$;} \\
0  , & \text{ if $\dim F(P) < d-1$,}
\end{cases}
$$
\end{thm}

\noindent
{\bf Proof.} The crepant  morphism  $\varphi\, :\, \widehat{Z} \to \widetilde{Z}$ implies the 
equalities  \[(K_{\widehat{Z}})^{d-1} = 
(\varphi^*(\widetilde{Z}))^{d-1} =  
(K_{\widetilde{Z}})^{d-1}. \]
Therefore it is sufficient to compute the top intersection number for the  canonical divisor on the 
canonical model $\widetilde{Z} \subset \widetilde{V}$. If 
the Kodaira dimension is not maximal, this number is $0$. 
The maximal Kodaira dimension $d-1$ 
of  $\widetilde{Z}$ appears 
only in two cases:  and $\dim F(P) = d-1$. 

Case 1.   $\dim F(P) = d$. In this case, by the adjunction formula,  we have 
\[ [K_{\widetilde{Z}}]^{d-1} = ([K_{\widetilde{V}}] + [\widetilde{Z}])^{d-1} \cdot [\widetilde{Z}] = 
([K_{\widetilde{V}}] + [\widetilde{Z}])^{d}  
- [K_{\widetilde{V}}] \cdot ([K_{\widetilde{V}}] + [\widetilde{Z}])^{d-1}.  \]
The intersection number $([K_{\widetilde{V}}] + [\widetilde{Z}])^{d}$ is the degree of the semiample adjoint 
$\Q$-divisor associated with the polytope $F(P)$. Therefore, 
we have 
\[   ([K_{\widetilde{V}}] + [\widetilde{Z}])^{d} = {\rm Vol}_d(F(P)). \]
In the computation of the intersection number 
$- [K_{\widetilde{V}}] \cdot ([K_{\widetilde{V}}] + [\widetilde{Z}])^{d-1}$
we use the formula
\[ - [K_{\widetilde{V}}] = \sum_{\nu \in 
\widetilde{\Sigma}[1]} [D_\nu] \]
and the fact that the degree of the restriction on a 
toric divisor $D_\nu \subset \widetilde{V}$  of the semiample $\Q$-divisor $[K_{\widetilde{V}}] + [\widetilde{Z}]$ associated with a polytope
$F(P)$ equals ${\rm Vol}_{d-1}(Q)$ if the supporting hyperplane with the normal vector $\nu$ defines a 
$(d-1)$-dimensional face of $F(P)$ and $0$ otherwise. 
It remains to take sum over all $\nu \in \widetilde{\Sigma}[1]$. 

Case 2.   $\dim F(P) = d-1$. In this case, the semiample 
adjoint $\Q$-divisor $[K_{\widetilde{V}}] + [\widetilde{Z}]$ on $\widetilde{V}$ associated with the polytope $F(P)$  defines a toric morphism $\vartheta\, :\, \widetilde{V} \to 
V_{F(P)}$ onto $(d-1)$-dimensional toric variety. Since 
the restriction of $\vartheta$ on the canonical model $\widetilde{Z}$ defines  a double covering $\widetilde{Z} \to 
V_{F(P)}$ we obtain 
\[ [K_{\widetilde{Z}}]^{d-1} = [\widetilde{Z}] \cdot 
 ([K_{\widetilde{V}}] + [\widetilde{Z}])^{d-1} = 2{\rm Vol}_{d-1}(F(P)). \]
\hfill $\Box$ 

\begin{ex}
Let $P$ be the $3$-dimensional simplex which is the Newton polytope of the equation of the 
Godeaux surface $S$ in $\P^3/\mu_5$ obtained as quotient 
of Fermat quintic surface $z_1^5 + z_2^5 + z_3^5 + z_4^5 =0$ modulo the action of the 5-order cyclic group $\mu_5 = 
\langle \zeta \rangle$:
\[ (z_1: z_2: z_3: z_4) \mapsto (\zeta z_1: \zeta^2 z_2: \zeta^3 z_3: \zeta^4 z_4). \]
Then the Fine interior 
$F(P)$ is a $3$-dimensional rational simplex with 
${\rm Vol}_3(F(P)) = 1/5$, and all its 4 facets $Q$ are rational 
triangles with  ${\rm Vol}_2(Q) = 1/5$. Thus, the  above formula yields $K_S^2 = 5 \times 1/5 = 1$. 
\end{ex}

\section{Further developments}

\subsection{Minimal surfaces} 
It is natural to study  
minimal surfaces arising from lattice polytopes $P$ of dimension $3$.  If $P \subset \R^3$
is a $3$-dimensional lattice polytope with $F(P) \neq \emptyset$,  then the Chern numbers $c_1^2$ and $c_2$ of minimal compactifications $\widehat{S}$ of nondegenerate 
surfaces $S \subset (\C^*)$ with the Newton 
polytope $P$ are completely determined by its Fine interior $F(P)$, because the number 
$c_1^2(\widehat{S}) = K^2_{\widehat{S}} $ 
is one of three integers 
\[ 0, \;\;  2{\rm Vol}_2(F(P)), \;\; {\rm Vol}_3(F(P)) + \sum_{Q \prec F(P) \atop \dim Q =2} {\rm Vol}_2(Q) \]
and 
\[ c_2(\widehat{S}) = 12 \chi({\mathcal O}_{\widehat{S}}) - c_1^2(\widehat{S}), \;\; \text{where $\chi({\mathcal O}_{\widehat{S}}) = 1 + |F(P) \cap \Z^3|. $} 
\]
The proposed in this paper general method for finding minimal models was applied to   
all   $674\,688$  
three-dimensional Newton 
polytopes with only one interior lattice point $0$. These polytopes 
were classified by Kasprzyk \cite{Kas10}. There exist $665\,599$ almost reflexive $3$-dimensional 
lattice polytopes $P$ including  $4\,319$ reflexive ones characterized by the condition $F(P) = 0$\footnote{A $d$-dimensional lattice polytope $P$ is called 
{\em almost reflexive} if $F(P) =\{0 \}$ and $C(P)$ is reflexive}. These polytopes define families of $K3$-surfaces. The remaining 
$9\,089$ polytopes give rise to minimal surfaces $\widehat{Z}$
of positive  Kodaira dimension $\kappa >0$  with $p_g =1$:  
elliptic surfaces $(\kappa =1)$, Todorov and Kanev surfaces ($\kappa =2$)  \cite{BKS19}.  

It would be interesting in general to investigate the geography of minimal surfaces $\widehat{S}$ arising  from arbitrary $3$-dimensional lattice polytopes $P$ with $2$ or more interior 
lattice points. 

\subsection{Minimal $3$-folds}
As for $3$-folds, the complete classification of $473\,800\,776$ four-dimensional 
reflexive polytopes obtained by Kreuzer and Skarke \cite{KS02} 
gives rise to a lot of topologically different examples of smooth $3$-dimensional Calabi-Yau varieties \cite{AGHJN15}. Note that  
the complete list of four-dimensional 
almost reflexive polytopes is still unknown  and it is expected 
to be huge. Therefore, it would be better to investigate some qualitative 
properties of the corresponding $3$-dimensional  minimal 
Calabi-Yau models, i.e. about possible isolated terminal $cDV$-singularities and their description in terms of  the corresponding almost reflexive $4$-dimensional Newton polytope $P \subset \R^4$. 
For example, the  $4$-dimensional almost reflexive lattice 
polytope $P \subset \R^4$: $x_i \geq -1\;
(1 \leq i \leq 4)$, $x_1 + x_2 + x_3 + x_4 \leq 1$,  $x_1 \leq 2$
provides simplest examples of nondegenerate $3$-dimensional hypersurfaces whose minimal Calabi-Yau models  
are not smooth. These minimal models 
 are $3$-dimensional Calabi-Yau quintics in $\P^4$ 
with an isolated terminal $cDV$-singularity 
analytically isomorphic to  $z_1^2 + z_2^2 + z_3^2 + z_4^2 =0$
\cite[Exam. 4.14]{Bat17} and \cite[Exam. 2.9]{BKS19}.

It was shown in \cite{BKS19} that up to unimodular equivalence there exist 
exactly $5$ three-dimensional Newton polytopes
defining Enriques surfaces. 
 It looks reasonable  to extend this classification in dimension $4$  and to obtain a complete list of 
all $4$-dimensional 
lattice polytopes $P$  with $\dim F(P) =0$ and 
$F(P) \cap M = \emptyset$. 
 This would yield  interesting examples of minimal 
$3$-folds with $\kappa =0$ that  are $3$-dimensional
analogs of Enriques surfaces (see \ref{exam-parall}). 
A particular  example of such a 
$4$-dimensional lattice polytope $P$
with  $F(P) = 
\{ (1/2,1/2,1/2,1/2) \}$ was proposed to author by Harald 
Skarke as  convex hull of the following 
$10$ lattice points in $\R^4$: 
\[ (0,0,0,0),\, (1,1,1,1),\, (2,0,0,0),\, (0,2,0,0),\, (0,0,2,0),\, (0,0,0,2), \] 
\[  (-1,1,1,1),\, (1,-1,1,1,),\, (1,1,-1,1),\, (1,1,1,-1). \]
The normal fan $\Sigma_P$ of $P$ is generated by $12$ primitive 
lattice vectors 
\[ \pm (1, 1, 0, 0), \; \pm (1,0,0,1), \; \pm (0,1,1,0), \; \pm (
0,1,0,1), \pm(1,0,1,0), \pm (0,0,1,1). \] 

\subsection{Mirror symmetry beyond reflexive polytopes}
The proposed method significantly extends our possibilities and 
allows to  construct minimal Calabi-Yau models 
of toric hypersurfaces  in case when  
the Newton polytope $P$ is not reflexive,  but 
only satisfies the condition 
$F(P) =0$ \cite{Bat17}.  Many examples of such lattice polytopes  
$P$ are contained   among  almost reflexive lattice 
polytopes $P$ of dimension $3$ and $4$ \cite[\S 2]{BKS19}. 
If $F(P) = 0$, then the canonical Calabi-Yau model $\widetilde{Z}$ is  the Zariski closure of 
$Z$ in the $\Q$-Gorenstein 
toric Fano variety $\widetilde{V}$ corresponding to the 
canonical hull  $C(P)$  of  $P$.  It is still an open problem to 
find a criterion for $d$-dimensional lattice 
polytope $P$ with $F(P) = \{0\}$ such that Calabi-Yau compactifications of nondegenerate affine toric 
hypersurfaces $Z \subset (\C^*)^d$ 
with the Newton polytope $P$ admit mirrors.

\subsection{Minimal models and homogeneous coordinates}
Let $P$ be a $d$-dimen\-sional Newton polytope with $F(P) \neq \emptyset$. We set $n := |S_F(P)|$. Then the 
simplicial terminal toric varieties $\widehat{V}$ defined above 
for constructing minimal models of nondegenerate hypersurfaces can be obtained 
as GIT-quotients of $\C^n$ by the 
Neron-Severi quasi-torus of dimension 
$n-d$. It is natural to investigate the proposed above combinatorial 
construction for 
canonical and minimal models of toric hypersurfaces 
using the homogeneous coordinates $\{ z_1,\ldots, z_n \}$ of 
$\widehat{V}$ \cite{Cox95}.  
This point of view possibly  allows to obtain  
an alternative interpretation of the Fine 
interior $F(P)$ via Newton polytopes of partial derivatives 
of the defining multi-homogeneous polynomial $h \in \C[z_1, \ldots, z_n]$ that defines the projective toric hypersurface 
$\{ h=0 \} \subset \widehat{V}$ using $n$ homogeneous coordinates \cite{BC94}. The homogeneous 
coordinates could allow to 
construct minimal models $\widehat{Z}$ under some weaker nondegeneracy  condition for coefficients of $h$ 
using the 
nonvanishing of the $A$-discriminant instead of  
the nonvanishing of the principal $A$-determinant  \cite{GKZ94}.  

\subsection{Complete intersections} 
The genus formula of Khovanski\v{i} for 
nondegenerate complete intersections \cite{Kho78} and 
the combinatorial construction of minimal models 
of Calabi-Yau complete intersections 
\cite{BB96a,BB96b,BL15} motivate natural 
generalizations of the combinatorial construction 
of canonical and minimal models 
of nondegenerate toric hypersurfaces 
to the case of complete intersections defined by $r$ Newton 
polytopes $P_1, \ldots, P_r \subset M_\R$. The crucial 
combinatorial condition $F(P) \neq \emptyset$ has to be 
applied to the Minkowski sum $P:= P_1 + \cdots + P_r$.
The canonical model of the affine 
complete intersection $Z = \bigcap_{i=1}^r Z_i \subset \T^d$  should be obtained 
as Zariski closure of $Z$ 
in the $\Q$-Gorenstein 
canonical toric variety $\widetilde{V}$ defined by the normal 
fan of the Minkowski sum 
\[ \widetilde{P} = C(P) + F(P) = C(P_1) + \cdots + C(P_r) + F(P), \]
where 
\[ C(P_i):= \{ x \in M_\R \, : \, \langle x, \nu \rangle\geq 
\ord_{P_i}(\nu) \;\; \; \forall \nu \in S_F(P)\}, \;\; i =1, \ldots, r. \]
This would provide  an alternative point of view on the 
work of Fletcher \cite{Fle00}.

\end{document}